\documentclass[hidelinks,onefignum,onetabnum]{siamart250211}

\usepackage{lipsum}
\usepackage{amsfonts}
\usepackage{graphicx}
\usepackage{epstopdf}
\usepackage{algorithmic}
\usepackage{subcaption}
\usepackage{tikz}
\usepackage{pgfplots}
\usepackage{siunitx}
\usetikzlibrary{3d,calc,arrows.meta,decorations.pathreplacing}
\usepackage{tikz-3dplot}
\pgfplotsset{compat=1.17}
\usepackage{xcolor}
\usepackage{subcaption} 
\usepackage{caption}
\usepackage{listings}
\usepackage{threeparttable}

\pgfplotsset{compat=1.18}
\ifpdf
  \DeclareGraphicsExtensions{.eps,.pdf,.png,.jpg}
\else
  \DeclareGraphicsExtensions{.eps}
\fi

\makeatletter
\AddToHook{cmd/appendix/before}{\def\cref@section@alias{appendix}}
\makeatother

\definecolor{lightcopper}{rgb}{0.85, 0.57, 0.35}
\definecolor{aluminum}{rgb}{0.65, 0.7, 0.75}
\definecolor{mutedblue}{RGB}{68,119,170}   
\definecolor{mutedorange}{RGB}{230,159,0}   
\definecolor{mutedgray}{RGB}{102,102,102}  

\newsiamremark{remark}{Remark}
\newsiamremark{hypothesis}{Hypothesis}
\crefname{hypothesis}{Hypothesis}{Hypotheses}
\newsiamthm{claim}{Claim}
\newsiamremark{fact}{Fact}
\crefname{fact}{Fact}{Facts}

\headers{Scalable Preconditioners for the P4D DFN Battery Model}{Thomas Roy et al.}

\title{{Scalable Preconditioners for the Pseudo-4D DFN Lithium-ion Battery Model}\thanks{Uploaded on February 6, 2026.
  \funding{
This work was performed under the auspices of the U.S. Department of Energy by Lawrence Livermore National Laboratory under Contract DE-AC52-07NA27344.
This work was supported by the Lawrence Livermore National Laboratory LDRD 23-SI-002.
LLNL release number: LLNL-JRNL-2015625-DRAFT.
}}}

\author{
  Thomas Roy\thanks{Lawrence Livermore National Laboratory, Livermore, CA, United States of America. Corresponding author: Thomas Roy (\email{roy27@llnl.gov}).}
  \and
  Nicholas W.~Brady\footnotemark[2]
  \and
  Giovanna Bucci\footnotemark[2]
  \and
  Nicholas R.~Cross\footnotemark[2]
  \and
  Victoria M.~Ehlinger\footnotemark[2]
  \and
  Tiras Y.~Lin\footnotemark[2]
  \and
  Hanyu Li\footnotemark[2]
  \and
  Marcus A.~Worsley\footnotemark[2]
}

\usepackage{amsopn}

\ifpdf
\hypersetup{
  pdftitle={Scalable solvers for the pseudo-4D battery model},
  pdfauthor={Thomas Roy}
}
\fi

\begin{document}

\maketitle

\begin{abstract}
The pseudo-4D Doyle–Fuller–Newman (DFN) model enables predictive simulation of lithium-ion batteries with three-dimensional electrode architectures and particle-scale diffusion, extending the standard pseudo-2D (P2D) formulation to fully resolve cell geometry.
This leads to large, nonlinear systems with strong coupling across multiple physical scales, posing significant challenges for scalable numerical solution.
We introduce block-structured preconditioning strategies that exploit the mathematical properties of the coupled system, employing multigrid techniques for electrode-level operators and localized solvers for particle-scale diffusion.
Comprehensive scalability studies are performed across a range of geometries, including homogeneous and heterogeneous cubic cells, flattened jelly-roll configurations, and triply periodic minimal surface electrodes, to assess solver robustness and parallel scalability.
The proposed methods consistently deliver efficient convergence and enable the solution of battery models with hundreds of millions of degrees of freedom on large-scale parallel hardware.
\end{abstract}



\section{Introduction}

Lithium-ion batteries are ubiquitous, with applications spanning consumer electronics, electric vehicles, and grid-scale energy storage systems.
Predictive, physics-based simulation is essential for the design, safety, and optimization of these devices \cite{brosa2022continuum}.
Among continuum-scale models for full battery cells, the Doyle--Fuller--Newman (DFN) model is a widely accepted standard \cite{doyle1993modeling,fuller1994simulation}.

The classical DFN model is often referred to as a pseudo-two-dimensional (P2D) formulation, as it couples one-dimensional through-cell transport with radial diffusion in spherical active-material particles, thereby capturing essential multiscale physics at manageable computational cost.
These assumptions are well suited to layered, homogeneous, slab-like electrodes, for which spatial variations are predominantly through the cell thickness and the geometry can be modeled effectively in one dimension.
However, many modern battery architectures increasingly violate this one-dimensional geometric assumption \cite{jung2014multi,mei2020three,cross2024viability}.

Commercial batteries are commonly manufactured as prismatic or cylindrical cells composed of thin electrodes that are stacked and rolled.
As a first approximation, P2D models are typically employed for such cells, neglecting variations along the length of the electrodes.
Nevertheless, tab placement and current-collector geometry can induce pronounced three-dimensional variations in electrochemical response \cite{samba2014impact,waldmann2023effects}.
In addition, certain forms of microstructural heterogeneity can be incorporated into continuum-scale simulations \cite{mistry2020stochasticity,parmananda2022probing}.
These conventional slab-like electrode geometries arise from traditional manufacturing techniques, while advanced manufacturing approaches enable new electrode architectures and increased thicknesses \cite{xue2026architecting,batista2023design,zhang20173d, xue2024interpenetrated,wiesner2023additive,wang2026ultra}.
Such structured electrodes naturally introduce nontrivial three-dimensional spatial variations in transport and reaction processes.
Moreover, additional physical mechanisms, including thermal transport and mechanical stress due to swelling, act across the full three-dimensional cell geometry and across electronically insulating layers \cite{liu2017multiphysics,kantharaj2019heat,lin2024shape,bucci2016formulation}.
Similarly, electrolyte motion induced by electrode deformation is an increasingly studied phenomenon that can significantly impact performance and degradation in commercial cells \cite{aiken2023tracking,bond2025operando}.

Several recent studies have extended the classical DFN equations to fully three-dimensional geometries while retaining radial particle-scale diffusion, yielding pseudo-four-dimensional (P4D) formulations that resolve geometric heterogeneity and multiphysics effects inaccessible to one-dimensional through-cell models.
At the cell scale, P4D electrochemical–thermal formulations for layered geometries have been implemented in commercial multiphysics software, including COMSOL Multiphysics\textregistered{} \cite{di2025understanding,lin2023novel,mei2020three,comsol62}.
To reduce the cost of resolving repeated layered structures, reduced-order P4D frameworks have been introduced in Abaqus\textregistered{}, preserving three-dimensional resolution \cite{hahn2023reduced,abaqus2024}.
At the electrode-architecture scale, P4D formulations have been employed to study explicitly resolved three-dimensional porous and additively manufactured geometries, including lattices, fins, and TPMS-based structures, using COMSOL and the open-source code Firedrake \cite{gunnarsson2023architectured,cross2024viability,FiredrakeUserManual}.
Beyond electro\-chemical--thermal coupling, strongly coupled P4D frameworks incorporating mechanics, heat transport, and porous electrolyte flow have been proposed in Abaqus\textregistered{} for fully three-dimensional cell geometries \cite{kulathu2024three}.
P4D formulations have also been combined with operando measurements to resolve depth-dependent lithiation heterogeneity in structured electrodes \cite{weddle2025depth}.
In contrast to this predominantly application-driven literature, Xu and Cao analyzed the numerical structure of the P4D DFN equations and established optimal finite element convergence rates in spatial dimensions greater than one \cite{xu2024optimal}.
Nearly all existing P4D studies rely on direct solvers, which fundamentally limits problem size and restricts physical fidelity to small-scale demonstrations with fewer than $10^5$ mesh elements.
To our knowledge, the only work employing a scalable solver is \cite{cross2024viability}, which represents prior work by the present authors.

The P4D model gives rise to tightly coupled nonlinear systems that can become extremely large.
Electrolyte concentration, ionic potential, electronic potential, and solid-phase concentration are strongly coupled through nonlinear charge-transfer reactions within the electrodes.
This strong coupling, combined with a wide range of characteristic time scales, limits the effectiveness of segregated solution strategies, which typically require prohibitively small time steps; fully coupled implicit approaches are therefore preferred.
This limitation was recognized early in the context of pseudo-2D models, where a direct solver, referred to as \texttt{BAND(J)}, exploited the banded structure of the resulting linear systems \cite{newman1968numerical,Newman1998}.
For pseudo-4D models, however, direct solvers do not scale with the size of the electrode-level system, which arises from the discretization of three-dimensional differential operators.
Scalable solution strategies for such systems therefore require iterative methods equipped with effective preconditioning strategies \cite{wathen2015preconditioning}.

A natural approach for solving the coupled nonlinear systems arising from fully implicit time discretization is Newton's method combined with a Krylov subspace method for the linearized systems.
The efficiency of this approach critically depends on the availability of robust and scalable preconditioners.
Naive application of incomplete factorization or monolithic multigrid methods typically fails or exhibits poor scalability.
Fortunately, the DFN system exhibits a rich block structure that can be exploited in the design of effective preconditioners.
Several preconditioning strategies have been proposed for battery simulations that neglect particle-scale diffusion, often in the context of explicitly resolving electrode microstructures \cite{wu2002newton,allen2021segregated,fang2019parallel}.

In this work, we present scalable preconditioners for the fully implicit solution of the pseudo-4D DFN model in three-dimensional cell geometries, and demonstrate both strong and weak scaling with respect to refinement of the three-dimensional electrode mesh.
We consider block preconditioners that exploit the intrinsic structure of the coupled system and the mathematical properties of its constituent blocks.
The electrode-level blocks correspond to elliptic or parabolic operators and are therefore well suited to multigrid methods \cite{brandt1977multi}.
In contrast, the particle-level block consists of many small, local systems that can be addressed efficiently using localized solution strategies.
Both block-diagonal and block-triangular preconditioning approaches are examined, allowing for a trade-off between computational cost and fidelity of inter-equation coupling.
These preconditioning strategies enable the solution of large-scale, fully coupled three-dimensional battery simulations with over $10^7$ mesh elements.

In \cref{sec:gov}, we describe the pseudo-4D DFN model for three-dimensional cell geometries.
In \cref{sec:discretization}, we describe the spatial and temporal discretizations, using finite elements for the electrode-level equations and finite differences in spherical coordinates for the particle-scale diffusion equations, with fully coupled implicit time-stepping.
In \cref{sec:solvers}, we analyze the resulting block system and introduce the proposed preconditioning strategies.
Finally, in \cref{sec:results}, we present four test cases and perform scalability studies for large-scale problems.

\section{Governing equations}
\label{sec:gov}
\begin{figure}[htb!]
    \centering
    \begin{tikzpicture}[scale=0.015, x={(1cm,0cm)}, y={(0.6cm,0.3cm)}, z={(0cm,1cm)}, every node/.style={font=\normalsize}]

  \def\depth{175}       
  \def\width{175}       
  \def\cathode{70}
  \def\separator{35}
  \def\anode{70}
  \def\height{\cathode+\separator+\anode}

  \coordinate (O) at (0,0,0);
  \coordinate (C) at (0,0,\cathode);
  \coordinate (S) at (0,0,\cathode+\separator);
  \coordinate (A) at (0,0,\height);

  \fill[mutedblue!75!gray] (O) -- ++(0,0,\cathode) -- ++(\depth,0,0) -- ++(0,0,-\cathode) -- cycle;
  \fill[mutedblue!70!gray!95!black] (\depth,0,0) -- ++(0,\width,0) -- ++(0,0,\cathode) -- ++(0,-\width,0) -- cycle;
  \node at (\depth/2, 0, \cathode/2) { 
  \begin{tabular}{c@{\hspace{0.0em}}c}
    \begin{tabular}{c}
      positive \\ 
      electrode
    \end{tabular} & 
    ($\Omega_p$)
  \end{tabular}
  };

  \fill[mutedgray!20] (C) -- ++(0,0,\separator) -- ++(\depth,0,0) -- ++(0,0,-\separator) -- cycle;
  \fill[mutedgray!25] (\depth,0,\cathode) -- ++(0,\width,0) -- ++(0,0,\separator) -- ++(0,-\width,0) -- cycle;
  \node at (\depth/2, 0, \cathode + \separator/2) {separator ($\Omega_m$)};

  \fill[mutedorange!85!white] (S) -- ++(0,0,\anode) -- ++(\depth,0,0) -- ++(0,0,-\anode) -- cycle;
  \fill[mutedorange!90!white] (\depth,0,\cathode+\separator) -- ++(0,\width,0) -- ++(0,0,\anode) -- ++(0,-\width,0) -- cycle;
  \fill[mutedorange!95!white] (A) -- ++(\depth,0,0) -- ++(0,\width,0) -- ++(-\depth,0,0) -- cycle; 
  \node at (\depth/2, 0, \cathode + \separator + \anode/2) {
  \begin{tabular}{c@{\hspace{0.0em}}c}
    \begin{tabular}{c}
      negative \\ 
      electrode
    \end{tabular} & 
    ($\Omega_n$)
  \end{tabular}
  };

  \pgfmathsetmacro{\zoomx}{\depth}
  \pgfmathsetmacro{\zoomy}{0.5*\depth}
  \pgfmathsetmacro{\zoomz}{\cathode + \separator + 0.5*\anode}

  \coordinate (arrowstart) at (\zoomx,\zoomy,\zoomz);
  \coordinate (arrowend) at (357,0,224);
  \coordinate (arrowend2) at (357,0,142);
  \draw[thick,gray,dashed] (arrowstart) -- (arrowend);
  \draw[thick,gray,dashed] (arrowstart) -- (arrowend2);

    \begin{scope}[reset cm, shift={(5.65,0,2.76)}, x={(1, 0, 0)}, y={(0, 0, 1)}, scale=0.1]  
    \fill[mutedorange!90!white] (0,0) circle[radius=6.5];
    \draw[thick] (0,0) circle[radius=6.5];
    \draw[thick,-{Stealth[length=5pt]}] (0,0) -- (6.5,6.5); 
    \draw[thick,{Stealth[length=3pt]}-{Stealth[length=3pt]}] (0,0) -- (0,-6.5); 
    \node at (2.2, -3.25) {\small $R_n$};
    \node at (7.2, 7.2) {$r$};
  \end{scope}

  \pgfmathsetmacro{\zoomx}{\depth}
  \pgfmathsetmacro{\zoomy}{0.5*\depth}
  \pgfmathsetmacro{\zoomz}{0.5*\cathode}
  \coordinate (arrowstart) at (\zoomx,\zoomy,\zoomz);
  \coordinate (arrowend) at (357,0,124);
  \coordinate (arrowend2) at (357,0,42);
  \draw[thick,gray,dashed] (arrowstart) -- (arrowend);
  \draw[thick,gray,dashed] (arrowstart) -- (arrowend2);

    \begin{scope}[reset cm, shift={(5.65,0,1.25)}, x={(1, 0, 0)}, y={(0, 0, 1)}, scale=0.1]  
    \fill[mutedblue!80!white] (0,0) circle[radius=6.5];
    \draw[thick] (0,0) circle[radius=6.5];
    \draw[thick,-{Stealth[length=5pt]}] (0,0) -- (6.5,6.5); 
    \draw[thick,{Stealth[length=3pt]}-{Stealth[length=3pt]}] (0,0) -- (0,-6.5); 
    \node at (2.2, -3.25) {\small $R_p$};
    \node at (7.2, 7.2) {$r$};
  \end{scope}

  \node at (0.2*\width,1.2*\height,1.4*\height) {Electrodes};
  \node at (1.4*\width,1.2*\height,1.4*\height) {Particles};
\end{tikzpicture}

    \caption{Schematic of the pseudo-4D DFN battery model with 3D electrodes and 1D spherical particles.}
    \label{fig:p4d}
\end{figure}
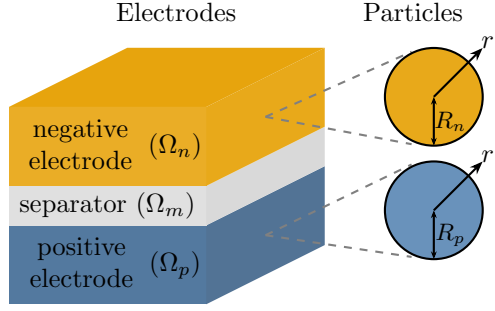

We consider the Doyle\textendash Fuller\textendash Newman (DFN) \linebreak model for lithium-ion batteries~\cite{doyle1993modeling,fuller1994simulation}, generalized to cells with three-dimensional geometries of arbitrary shape.
The domain $\Omega \subset \mathbb{R}^3$ is partitioned into negative electrode (anode) $\Omega_n$, separator (membrane) $\Omega_m$, and positive electrode (cathode) $\Omega_p$, such that $\Omega = \Omega_n \cup \Omega_m \cup \Omega_p$ and these subdomains are disjoint.
By convention, electrode signs in batteries are defined according to their polarity during discharge, when the cell delivers power.
Each point $\mathbf{x} \in \Omega_n \cup \Omega_p$ is associated with a representative spherical active material particle of radius $R_s(\mathbf{x})$.
Here we assume $R_s=R_p$ in $\Omega_p$ and $R_s=R_n$ in $\Omega_n$, where both $R_p$ and $R_n$ are constant values associated with the active material in their respective electrode.
This corresponds to a homogenized description in which particle size distributions are approximated by representative mean values.
The electrode-scale and particle-scale domains are coupled through charge-transfer reactions occurring throughout the electrodes and at particle surfaces, as detailed below, yielding the pseudo-4D model illustrated in \cref{fig:p4d}.
Parameter values for the following equations are given in \cref{sec:model_parameters} for the test cases considered in \cref{sec:results}.

\subsection{Lithium diffusion in active particles}

The solid-phase lithium concentration $c_s = c_s(r, \mathbf{x}, t)$ in the particle satisfies
\begin{equation} \label{eq:particle}
  \frac{\partial c_s}{\partial t} = \frac{1}{r^2} \frac{\partial}{\partial r} \left( D_s r^2 \frac{\partial c_s}{\partial r} \right), \quad r \in [0, R_s(\mathbf{x})],
\end{equation}
where \( D_s \) is the solid-phase diffusivity.
At the particle center, symmetry requires
\begin{equation} \label{eq:particleBC1}
  \left. \frac{\partial c_s}{\partial r} \right|_{r=0} = 0.
\end{equation}
At the surface, the radial flux is proportional to the reaction current density \( i_n \),
\begin{equation}\label{eq:particleBC2}
  -D_s \left. \frac{\partial c_s}{\partial r} \right|_{r=R_s} = \frac{i_n}{F},
\end{equation}
where $F$ is Faraday's constant.

\subsection{Electrolyte transport}

The electrolyte-phase lithium concentration \( c_e = c_e(\mathbf{x}, t) \) satisfies the conservation equation
\begin{equation}
  \varepsilon \frac{\partial c_e}{\partial t}
  - \nabla \cdot \left( D_e^{\mathrm{eff}} \nabla c_e \right)
  = \frac{1 - t_+^0}{F} a i_n, \quad \text{in } \Omega,
\end{equation}
where \( \varepsilon(\mathbf{x}) \) is the local porosity, \( t_+^0 \) is the transference number, and \(a (\mathbf{x})\) is the electrochemically active surface area per volume ratio.
The transference number may in general depend on the electrolyte concentration $c_e$ \cite{doyle1993modeling}; here it is taken to be constant, consistent with the material parameters specified in \cref{sec:model_parameters}.
The effective diffusivity is given by the Bruggeman correlation \cite{bruggeman1935berechnung}
\begin{equation}
  D_e^{\mathrm{eff}}(\mathbf{x}, c_e) = \varepsilon(\mathbf{x})^b D_e(c_e),
\end{equation}
where \( D_e(c_e) \) is the bulk electrolyte diffusivity and \( b = 1.5\) is the Bruggeman exponent.

\subsection{Charge conservation}

The electrolyte-phase (ionic) potential,\linebreak \( \phi_e = \phi_e(\mathbf{x}, t) \) satisfies
\begin{equation}
  -\nabla \cdot \left( \kappa^{\mathrm{eff}} \nabla \phi_e
  + \kappa^{\mathrm{eff}} \frac{2RT(1 - t_+^0)}{F c_e} \nabla c_e \right)
  = a i_n, \quad \text{in }\Omega.
\end{equation}
The effective electrolyte conductivity is defined as
\begin{equation}
  \kappa^{\mathrm{eff}}(\mathbf{x}, c_e) = \varepsilon(\mathbf{x})^b \kappa(c_e),
\end{equation}
where \( \kappa(c_e) \) is the bulk electrolyte conductivity.

The solid-phase (electronic) potential \( \phi_s = \phi_s(\mathbf{x}, t) \) is defined in the electrode domains and satisfies
\begin{equation}
  \label{eq:phis}
  -\nabla \cdot \left( \sigma^{\mathrm{eff}} \nabla \phi_s \right)
  = - a i_n, \quad  \text{in } \Omega_n \cup \Omega_p.
\end{equation}
The effective solid-phase conductivity accounts for the active material phase volume fraction \( \varepsilon_s (\mathbf{x}) \), and is given by
\begin{equation}
  \sigma^{\mathrm{eff}}(\mathbf{x}) = \varepsilon_s (\mathbf{x})^b \sigma,
\end{equation}
where \( \sigma \) is the bulk conductivity of the active material.

\subsection{Interfacial kinetics}

The interfacial reaction current density \( i_n  \) follows the Butler--Volmer relation in exponential form:
\begin{equation}
  i_n = i_0 \left[ \exp\left( \frac{\alpha_a F}{RT} \eta \right)
  - \exp\left( -\frac{\alpha_c F}{RT} \eta \right) \right],
\end{equation}
where \( \eta = \phi_s - \phi_e - U_{\mathrm{ocp}}(c_s^{\mathrm{surf}}) \) is the overpotential, \( R \) is the universal gas constant, and \( T \) is the temperature.
Here, \( \alpha_a \) and \( \alpha_c \) are the anodic and cathodic charge-transfer coefficients (\( \alpha_a + \alpha_c = 1 \)), and \( c_s^{\mathrm{surf}} = c_s(R_s(\mathbf{x}), \mathbf{x}, t) \) is the surface concentration of the active particle.
The exchange current density \( i_0 \) depends on both the electrolyte and solid-phase concentrations, typically modeled as
\begin{equation}
  i_0 = k \, c_e^{\alpha_a} \left(c_{s,\max} - c_s^{\mathrm{surf}}\right)^{\alpha_a}
  \left(c_s^{\mathrm{surf}}\right)^{\alpha_c},
\end{equation}
where \( k \) is a rate constant and \( c_{s,\max} \) is the maximum lithium concentration in the solid.

\subsection{Boundary and initial conditions}

Boundary conditions are imposed on \( \partial \Omega \) to model current collectors, insulation, and symmetry, depending on the physical configuration.
In this work, we will consider applied current boundary conditions (as opposed to applied voltage).
Thus, at the positive current collector, a Neumann condition is applied to \( \phi_s \) to enforce the applied current density such that
\begin{equation}
  \label{eq:Neumann}
  \sigma^{\mathrm{eff}}\nabla \phi_s \cdot \mathbf{n} = \frac{I_{\mathrm{app}}}{A},\quad \text{on } \Gamma_p
\end{equation}
where \( I_{\mathrm{app}} \) is the applied current and \( A \) is the area of the positive current collector \( \Gamma_p \).
At the negative current collector \( \Gamma_n \), the solid potential is grounded via a Dirichlet boundary condition:
\begin{equation}
  \label{eq:Dirichlet}
  \phi_s = 0, \quad \text{on } \Gamma_n.
\end{equation}
At the interfaces between the electrodes (\(\Omega_n\), \(\Omega_p\)) and the separator (\(\Omega_m\)), a zero-flux (Neumann) boundary condition is imposed on $\phi_s$.
The electrolyte phase (\(\phi_e\), \(c_e\)) is subject to zero-flux conditions on insulating boundaries (the exterior boundary of $\Omega$).
Initial conditions must be provided for \( c_e(\mathbf{x}, 0) \), \( \phi_e(\mathbf{x}, 0) \), \( \phi_s(\mathbf{x}, 0) \), and \( c_s(r, \mathbf{x}, 0) \).

\section{Discretization}
\label{sec:discretization}

In this section, we introduce our discretization scheme including a finite difference scheme for the particle equation and a finite element discretization for the electrode-level equations.
We also note details related to implementing the coupled scheme in a finite element framework.

\subsection{Spatial discretization of the particle equation}
\label{sec:particle-discretization}

We discretize the radial diffusion equation \cref{eq:particle} for the solid-phase concentration, \( c_s(r) \), in spherical active material particles using a finite difference method on a one-dimensional, nonuniform spherical mesh.
The scheme combines the scheme for spherical coordinates from \cite{thibault1987finite} and the scheme for nonuniform meshes from \cite{sundqvist1970simple}.

Let \( \{ r_i \}_{i=1}^{N_c} \) denote a strictly increasing sequence of radial node locations such that \( r_1 = 0 \) and \( r_{N_c} = R_s \), defining the computational mesh.
We define \( \Delta r_{i-1} = r_i - r_{i-1} \), \( \Delta r_i = r_{i+1} - r_i \), and the local mesh ratio \( \theta_i = \Delta r_i / \Delta r_{i-1} \).

The differential operator in \cref{eq:particle} is written as a sum of two terms:
\begin{equation}
\frac{1}{r^2} \frac{d}{dr} \left( r^2 \frac{d c_s}{dr} \right) = \frac{2}{r} \frac{d c_s}{dr} + \frac{d^2 c_s}{dr^2},
\end{equation}
each of which is approximated separately.

For interior nodes \( 2 \leq i \leq N_c - 1 \), we discretize the first derivative and second derivative terms as
\begin{equation}
\label{eq:L1}
\mathcal{L}_1(c)_i = \frac{2 D_s}{r_i}  \frac{c_{i+1} - \theta_i^2 c_{i-1} - (1 - \theta_i^2) c_i}{\Delta r_i (1 + \theta_i)},
\end{equation}
\begin{equation}
\label{eq:L2}
\mathcal{L}_2(c)_i = D_s \frac{2}{\Delta r_i \Delta r_{i-1}(1 + \theta_i)} \left( c_{i+1} + \theta_i c_{i-1} - (1 + \theta_i) c_i \right).
\end{equation}
The full discrete Laplacian at node \( i \) is given by
\begin{equation}
  \label{eq:Lh}
\mathcal{L}_h(c)_i = \mathcal{L}_1(c)_i + \mathcal{L}_2(c)_i.
\end{equation}

The zero-flux boundary condition at the particle center \cref{eq:particleBC1} is discretized using a second-order one-sided approximation:
\begin{equation}
  \label{eq:L0}
\mathcal{L}_h(c)_1 = \frac{6 D_s}{\Delta r_1^2} \left( c_2 - c_1 \right).
\end{equation}

At the particle surface, the flux boundary condition \cref{eq:particleBC2} is enforced using a second-order backward difference:
\begin{equation}
  \label{eq:LR}
\mathcal{L}_h(c)_{N_c} = \frac{2 D_s}{\Delta r_{N_c - 1}^2} \left( c_{N_c - 1} - c_{N_c} \right) + \frac{2 i_n}{F \Delta r_{N_c - 1}} \left( 1 + \frac{\Delta r_{N_c - 1}}{R_s} \right).
\end{equation}

The spatially discretized particle system is thus given by

\begin{equation} \label{eq:FD}
  \frac{\partial c_{s, i}}{\partial t} = \mathcal{L}_h (c_s)_i, \quad \text{for } i = 1,\dots,N_c,
\end{equation}
where \( c_{s, i} \approx c_s (r_i) \).

\begin{remark}
\label{rem:reg}
From \cite{sundqvist1970simple}, the spatial discretization described above is second-order accurate in space, provided the mesh is locally regular.
Specifically, if adjacent cells satisfy \(\Delta r_i - \Delta r_{i-1} = \mathcal{O}(\Delta r_{i-1}^2) \) or \(\mathcal{O}(\Delta r_{i}^2) \), then the local truncation error is \( \mathcal{O}(\Delta r_{i-1}^2) \) at all nodes.
If this condition is violated (e.g., under strongly graded meshes), the scheme remains consistent but only achieves first-order spatial accuracy.
\end{remark}

For this study, we choose a regular mesh defined using a geometric progression.
Let the ratio between the first and last element \(\bar \theta = \Delta r_1 / \Delta r_{N_c-1} \), where \(0 < \bar \theta < 1\).
The elements get progressively smaller using the relationship
\begin{equation}
  \Delta r_i = \bar{\theta}^{\frac{1}{N_c-2}} \Delta r_{i-1}.
\end{equation}
This exponentially smooth mesh is more regular than the regularity required in Remark \ref{rem:reg}.

\subsection{Spatial discretization for the electrode-level equations}
\label{sec:fem-electrode}

The three-dimensional equations for the electrolyte concentration \( c_e \), electrolyte potential \( \phi_e \), and solid potential \( \phi_s \) are discretized using continuous, piecewise linear finite elements on a conforming tetrahedral or hexahedral mesh of the domain \( \Omega \).
Let \( V_h \subset H^1(\Omega) \) denote the standard finite element space of continuous, piecewise linear functions, and \( V_h^0=\left\{ v\in V_h \mid v = 0 \text{ on }\Gamma_n \right\}\), i.e. to strongly enforce the ground boundary condition \cref{eq:Dirichlet}.

The governing equations are coupled through the reaction current density \( i_n \), which depends on the solution of the particle-level diffusion equation at each point in the electrodes.
We adopt a monolithic approach in which all variables including the particle concentrations are solved simultaneously as part of a single nonlinear system.
Here we write the variational formulation for the electrode-level equations.

We seek \( c_e, \phi_e \in V_h \) and \( \phi_s \in V_h^0 \) such that the following variational equations are satisfied for all test functions \( v_{c_e}, v_{\phi_e} \in V_h \) and \( v_{\phi_s} \in V_h^0 \), respectively:

\begin{equation}
\label{eq:ce-weak}
\int_{\Omega} \varepsilon \frac{\partial c_e}{\partial t} v_{c_e} \, \mathrm{d}\mathbf{x}
+ \int_{\Omega} D_e^{\mathrm{eff}} \nabla c_e \cdot \nabla v_{c_e} \, \mathrm{d}\mathbf{x}
= \int_{\Omega} \frac{1 - t_+^0}{F} a i_n \, v_{c_e} \, \mathrm{d}\mathbf{x},
\end{equation}
\begin{equation}
\label{eq:phie-weak}
\int_{\Omega} \kappa^{\mathrm{eff}} \nabla \phi_e \cdot \nabla v_{\phi_e} \, \mathrm{d}\mathbf{x}
+ \int_{\Omega} \kappa^{\mathrm{eff}} \frac{2RT(1 - t_+^0)}{F c_e} \nabla c_e \cdot \nabla v_{\phi_e} \, \mathrm{d}\mathbf{x}
= \int_{\Omega} a i_n \, v_{\phi_e} \, \mathrm{d}\mathbf{x},
\end{equation}
\begin{equation}
\label{eq:phis-weak}
\int_{\Omega_n \cup \Omega_p} \sigma^{\mathrm{eff}} \nabla \phi_s \cdot \nabla v_{\phi_s} \, \mathrm{d}\mathbf{x}
= -\int_{\Omega_n \cup \Omega_p} a i_n \, v_{\phi_s} \, \mathrm{d}\mathbf{x}
+ \int_{\Gamma_p} \frac{I_{\mathrm{app}}}{A}\, v_{\phi_s} \, \mathrm{d}\mathbf{s},
\end{equation}
while simultaneously solving for \( c_s \) and satisfying \cref{eq:FD}.

\subsection{Temporal discretization}

For time-stepping, we use backward Euler, \linebreak which is fully implicit and first-order accurate.

The accumulation terms for the liquid-phase and solid-phase concentrations are approximated by
\begin{equation}
  \frac{\partial c_s}{\partial t} \approx \frac{c_s^{i+1} - c_s^{i}}{t_{i+1}-t_i},
\end{equation}
\begin{equation}
  \frac{\partial c_e}{\partial t} \approx \frac{c_e^{i+1} - c_e^{i}}{t_{i+1}-t_i},
\end{equation}
where \( c_s^i \), \( c_e^i \) are the approximate solutions at time $t_i$.
The remaining terms in the spatially discretized equations are evaluated at time \( t_{i+1} \) so that the scheme is implicit.

While the spatial discretization described above is relatively simple, it is likely sufficient in most cases.
In contrast, for the temporal discretization, substantial improvements in numerical accuracy and computational efficiency could be achieved through the use of higher-order and adaptive time-stepping schemes, respectively.
For example, adaptive high order schemes from Sundials \cite{hindmarsh2005sundials,gardner2022sundials} are used in PyBaMM \cite{Sulzer2021}, leading to very fast and stable pseudo-2D simulations.
Since we are focusing on scalable solvers for pseudo-4D simulations, we only consider backward Euler for simplicity.

\subsection{Implementation details}
\label{sec:implementation}

The discretized equations are implemented using Firedrake \cite{FiredrakeUserManual}, an open-source finite element framework that enables the specification of variational formulations at a high level of abstraction.
Firedrake integrates natively with PETSc \cite{petsc-efficient,petsc-user-ref,petsc-web-page}, providing access to scalable nonlinear and linear solvers, including the solver strategies discussed in \cref{sec:solvers}.
An additional motivation for this choice is the availability of automated adjoint evaluation through pyadjoint \cite{mitusch2019dolfinadjoint}, which is required for PDE-constrained design optimization.
Firedrake has previously been used to perform topology optimization for the design of porous electrodes for electrochemical devices \cite{roy2022topology,li2024topology,barzegari2025topology}.
The code developed for this work is also used in ongoing studies on topology-optimizing battery electrodes.

At the time this work was conducted, Firedrake did not robustly support finite element spaces restricted to subdomains.
As a result, all variables are defined over the full computational domain.
In particular, the solid-phase potential and solid-phase concentration are defined everywhere, even though they are physically relevant only within electrode regions.
This design choice is necessary for topology optimization, where electrode locations evolve during the optimization process \cite{roy2022topology,li2024topology}.
In contrast, many battery model implementations restrict these variables to fixed positive and negative electrode subdomains.

To avoid matrix singularities in separator or electrolyte-only regions ($\Omega_m$), a small regularization value of \(\SI{1e-22}{S/m}\) is assigned to the effective solid-phase conductivity \( \sigma^\mathrm{eff} \) in these regions, whereas physical conductivity values for the electrodes are on the order $ > \SI{1}{S/m}$.
This also effectively mimics the no-flux boundary condition at the interface between $\Omega_m$ and the electrode domains ($\Omega_n$, $\Omega_p$).
Similarly, the particle diffusion equation is defined such that \( \partial c_s / \partial t = 0 \) in regions outside the electrodes.
These regularization choices were verified to have no measurable impact on the solution within the electrode regions.
Recent developments in Firedrake now allow finite element spaces to be defined on subdomains \cite{sagiyama2025abstraction}, and future implementations may exploit this capability when topology optimization is not performed.

The solid-phase concentration diffusion equation is incorporated into the finite element framework via a variational formulation over the particle degrees of freedom; details are provided in \cref{sec:appendix}.
Since standard continuous finite element spaces would enforce artificial continuity across the electrode domain, a piecewise constant finite element space is used for the solid-phase concentration.
Specifically, each solid-phase concentration component corresponding to a discrete particle radius \( r = r_i \) is represented as a cellwise constant field on the electrode-scale mesh.
While this choice leads to numerous particle-level unknowns on unstructured tetrahedral meshes, it yields a natural block structure that is well suited to the solver strategies described in \cref{sec:solvers}.

\section{Solution strategy}
\label{sec:solvers}

In large-scale simulations involving three-dimensional PDEs, the resulting linear systems can be too large to be efficiently solved with direct solvers due to their prohibitive memory requirements and computational cost, especially at fine spatial resolutions.
Due to the complexity of the P4D model, the limits of direct solvers are reached quickly.
Iterative solvers offer a scalable alternative, enabling the solution of large, sparse linear systems with more favorable memory usage and better suitability for parallel computing environments.
This section outlines the structure of the linear systems arising from the P4D model discretization and details the preconditioning strategies employed to achieve scalable simulations.

\subsection{Linear system structure}
After discretization, we linearize the nonlinear system \cref{eq:phis-weak,eq:phie-weak,eq:ce-weak,eq:FD} using Newton's method.
The nonlinear residual can be written as the block residual equation
\begin{equation}
  \mathcal{R} = \begin{bmatrix} R_{\phi_s}\\
                          R_{\phi_e}\\
                          R_{c_e}\\
                          R_{c_s}
  \end{bmatrix}
              = 0,
\end{equation}
where \( R_{\phi_s} \) is the residual of \cref{eq:phis-weak}, \( R_{\phi_e} \) is the residual of \cref{eq:phie-weak}, \( R_{c_e} \) is residual of \cref{eq:ce-weak}, and \( R_{c_s} \) is the residual of \cref{eq:FD}.

The linearization of this system results in a block matrix of this form:
\begin{equation}
  \label{eq:blocks}
\begin{bmatrix}
A_{\phi_s,\phi_s} & A_{\phi_s,\phi_e} & A_{\phi_s,c_e} & A_{\phi_s,c_s} \\
A_{\phi_e,\phi_s} & A_{\phi_e,\phi_e} & A_{\phi_e,c_e} & A_{\phi_e,c_s} \\
A_{c_e,\phi_s}    & A_{c_e,\phi_e}    & A_{c_e,c_e}    & A_{c_e,c_s}    \\
A_{c_s,\phi_s}    & A_{c_s,\phi_e}    & A_{c_s,c_e}    & A_{c_s,c_s}
\end{bmatrix}
,
\end{equation}
where each block \( A_{ij} \) represents the linearization of residual equation \( R_i \) with respect to variable \( j \in [\phi_s, \phi_e, c_e, c_s] \).
The blocks related to \( c_s \) are ordered point-wise according to the electrode dimension.
For example, in the case of the diagonal block
\begin{equation}
    A_{c_s, c_s} = \begin{bmatrix}
                    (A_{c_s, c_s})_{1,1} & \cdots & (A_{c_s, c_s})_{n,n} \\
                  \vdots & \ddots & \vdots \\
                  (A_{c_s, c_s})_{n,1} & \cdots & (A_{c_s, c_s})_{n,n}
                 \end{bmatrix}
                =
                \begin{bmatrix}
                    (A_{c_s, c_s})_{1,1} &  &  \\
                   & \ddots &  \\
                   & & (A_{c_s, c_s})_{n,n}
                 \end{bmatrix},
\end{equation}
where \( (.)_{i,j} \) represents the linearization at the element \( i \) with respect to element \( j \).
Since there is no dependence between electrode-level elements for the particle equation, this is a block diagonal matrix.
Each sub-block \( (A_{c_s, c_s})_{i,i} \) is an \( N_c \times N_c \) tri-diagonal matrix for the 1D finite-difference particle system at element \( i \).

\begin{figure}[htbp]
  \centering
  \includegraphics[width=0.8\linewidth]{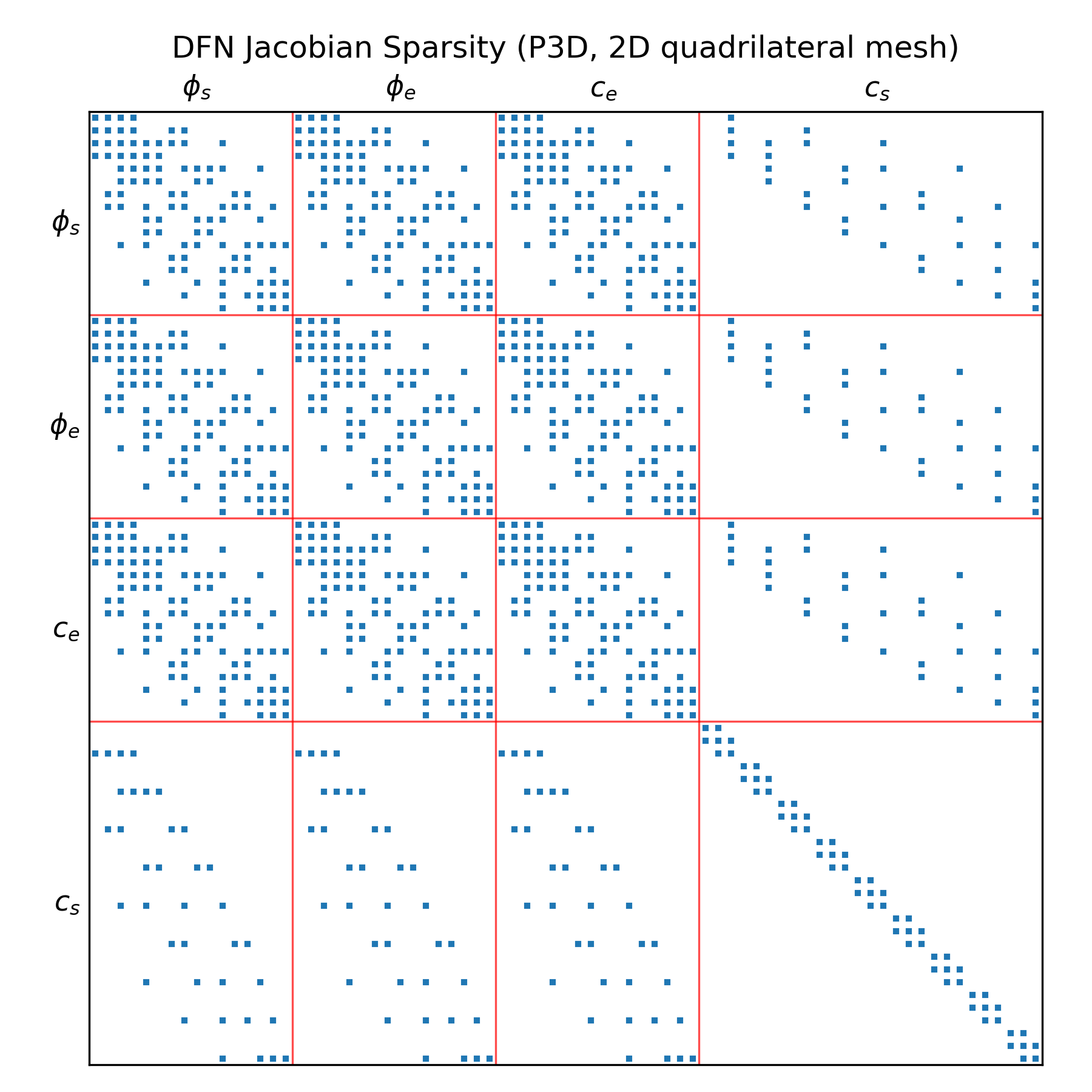}
  \caption{Sparsity pattern of the linearized system matrix for a P3D case with a \( 3 \times 3 \) quadrilateral mesh for the electrode-level domain and 3 points in the particle direction.}
  \label{fig:sparsity}
\end{figure}

We illustrate in \cref{fig:sparsity} the sparsity of the full block matrix from \cref{eq:blocks} for a pseudo-3D (P3D) case with a \( 3 \times 3 \) quadrilateral mesh for the electrode-level domain and 3 points in the particle direction.
We easily see the more connected structure of the finite element discretization for the top left blocks (the connectivity is even higher in 3D).
The bottom right block consists of small tri-diagonal matrices, one for each element in the electrode-level mesh.
The bottom left and top right sections correspond to the couplings between the electrode-level equations and the particle equation; they are sparse because they are only coupled through the solid-phase concentration at the surface.

\subsection{Preconditioning}
\label{sec:preconditioning}

The linearized systems are solved iteratively using GMRES~\cite{saad1986gmres}.
To ensure scalability of the solution procedure, we will employ appropriate preconditioning strategies~\cite{wathen2015preconditioning}.
A good preconditioner is typically an approximation of the system matrix for which the inverse can be efficiently computed or approximated.
For systems of coupled PDEs, considering the block structure is often key to designing a preconditioner since most available solvers will not work well when applied to the monolithic system.

A simple strategy for block systems is to consider a block diagonal preconditioner, that is, to ignore all coupling between equations:
\begin{equation}
  \label{eq:PJ}
P_{J} =
\begin{bmatrix}
A_{\phi_s,\phi_s} & 0                   & 0               & 0 \\
0                 & A_{\phi_e,\phi_e}   & 0               & 0 \\
0                 & 0                   & A_{c_e,c_e}     & 0 \\
0                 & 0                   & 0               & A_{c_s,c_s}
\end{bmatrix}.
\end{equation}
Inverting \( P \) only requires applying the inverse of each of its blocks:
\begin{equation}
P_{J}^{-1} =
\begin{bmatrix}
A_{\phi_s,\phi_s}^{-1} & 0                        & 0                    & 0 \\
0                      & A_{\phi_e,\phi_e}^{-1}   & 0                    & 0 \\
0                      & 0                        & A_{c_e,c_e}^{-1}     & 0 \\
0                      & 0                        & 0                    & A_{c_s,c_s}^{-1}
\end{bmatrix}.
\end{equation}
This preconditioner is often referred to as block Jacobi (BJ), or in PETSc as \texttt{additive fieldsplit}.

A more expensive option is to only ignore half the couplings, leading to a block triangular preconditioner:
\begin{equation}
  \label{eq:PGS}
P_{GS} =
\begin{bmatrix}
A_{\phi_s,\phi_s} & 0                     & 0                  & 0 \\
A_{\phi_e,\phi_s} & A_{\phi_e,\phi_e}     & 0                  & 0 \\
A_{c_e,\phi_s}    & A_{c_e,\phi_e}        & A_{c_e,c_e}        & 0 \\
A_{c_s,\phi_s}    & A_{c_s,\phi_e}        & A_{c_s,c_e}        & A_{c_s,c_s}
\end{bmatrix}.
\end{equation}
While inverting \( P_{GS} \) requires more operations, it also only needs the inverse of each diagonal block:
\begin{equation}
P_{GS}^{-1} =
\begin{bmatrix}
A_{\phi_s,\phi_s}^{-1} & 0 & 0 & 0 \\
- A_{\phi_e,\phi_e}^{-1} A_{\phi_e,\phi_s} A_{\phi_s,\phi_s}^{-1} & A_{\phi_e,\phi_e}^{-1} & 0 & 0 \\
\ast & \ast & A_{c_e,c_e}^{-1} & 0 \\
\ast & \ast & \ast & A_{c_s,c_s}^{-1}
\end{bmatrix},
\end{equation}
where the \(\ast\) entries represent nested expressions involving combinations of previously inverted blocks.
For example:
\begin{equation}
(P_{GS}^{-1})_{3,1} = -A_{c_e, c_e}^{-1} \left( 
    A_{c_e, \phi_s} A_{\phi_s, \phi_s}^{-1} 
    + A_{c_e, \phi_e} A_{\phi_e, \phi_e}^{-1} A_{\phi_e, \phi_s} A_{\phi_s, \phi_s}^{-1} 
\right),
\end{equation}
and so on.
This method is often referred to as block Gauss-Seidel (BGS), or in PETSc as \texttt{multiplicative fieldsplit}.

Note that for a \( 4 \times 4 \) block matrix, there are 24 different choices of ordering and the ordering shown above is arbitrary.
In \cref{sec:ordering}, we will briefly discuss the impact of system ordering on the performance of \( P_{GS} \).

To achieve a scalable method, the inverse of the diagonal blocks must be approximated using an iterative method instead of directly inverting them.
We note that the potential blocks, \( A_{\phi_s, \phi_s} \) and \( A_{\phi_e, \phi_e} \), are elliptic operators, and that the electrolyte concentration block \( A_{c_e, c_e} \) is a parabolic operator, meaning they are ideal candidates for a multigrid method \cite{brandt1977multi}.
As for the solid-phase concentration block \( A_{c_s, c_s} \), it consists of small \( N_c \times N_c \) tri-diagonal matrices, corresponding to each mesh element.
The inverse can be computed by simply inverting each small matrix, which can be done efficiently with a direct solver because these are small 1D systems.
Since each small matrix corresponds to a system defined on a single element, this can be done using element-wise block Jacobi, or point block Jacobi in PETSc.
Given that \( A_{c_s, c_s}, A_{c_s, c_s} \) is itself tri-diagonal, symmetric Gauss-Seidel is also a viable option as a preconditioner (results not shown here).

In \cref{sec:results}, we will consider the scalability of the BJ and BGS preconditioners defined in \cref{eq:PJ} and \cref{eq:PGS}, respectively.
In both cases the inverses of the electrode-level blocks are approximated with single Algebraic Multigrid (AMG) V-cycles \cite{ruge1987algebraic}:
\begin{equation}
  A_{\phi_s, \phi_s}^{-1} \approx \mathrm{AMG}(A_{\phi_s,\phi_s}), \quad
  A_{\phi_e, \phi_e}^{-1} \approx \mathrm{AMG}(A_{\phi_e,\phi_e}), \quad
  A_{c_e, c_e}^{-1} \approx \mathrm{AMG}(A_{c_e,c_e}), 
\end{equation}
and the inverse of the particle block is computed using element-wise block Jacobi (EBJ):
\begin{equation}
  A_{c_s, c_s}^{-1} = \mathrm{EBJ}(A_{c_s,c_s}).
\end{equation}

For AMG, we use BoomerAMG~\cite{henson2002boomeramg} from the hypre library~\cite{hypre}.
The full PETSc solver parameters are given in \cref{sec:solver_parameters}.

As noted in \cref{sec:implementation}, the solid-phase potential $\phi_s$ and the solid-phase concentration $c_s$ are defined and solved over the whole domain (even the separator).
If these are instead defined as separate variables local to each cathode, one would get four electrode-specific blocks instead of $A_{\phi_s,\phi_s}$ and $A_{c_s,c_s}$.
The preconditioning strategies detailed here could easily be adapted to this formulation.
The additive approach with $P_{J}$ would essentially be the same since the variables defined in one electrode are not directly coupled with the ones in the other electrode.
The multiplicative approach with $P_{GS}$ would be slightly different if electrode-specific equations are solved sequentially.
In the case of simulating resolved microstructures, block preconditioners with electrode-specific blocks are considered in \cite{fang2019parallel}.


\section{Results}
\label{sec:results}
In this section, we evaluate block Jacobi (BJ) and block \linebreak Gauss\textendash Seidel (BGS) preconditioners for Pseudo-4D DFN simulations.
We first describe four test cases with increasing spatial complexity.
We then analyze block ordering effects for BGS.
Finally, we present weak and strong scaling results for both preconditioners.

\subsection{Test cases}
\label{sec:cases}

We consider four different test cases with varying levels of complexity in terms of 3D spatial variations: a cubic cell with homogeneous properties (Case I), then with heterogeneous properties (Case II), a flattened jelly roll leading to a very skewed geometry (Case III), and interpenetrating gyroid electrodes with a complex geometry (Case IV).

For simplicity, all cases use the material properties from \cite{marquis2019asymptotic}, corresponding to a lithium cobalt oxide cathode (positive electrode) and a graphite anode (negative electrode).
All model parameters are detailed in \cref{sec:model_parameters}.

All simulations are performed for a fixed applied current density corresponding to a battery discharge at a rate of 1C.
This corresponds to a theoretical full discharge in one hour, although a minimum voltage (for safety) would be reached much sooner than that.
For simplicity, the simulated time is 30 minutes, which means half of a full discharge.
This avoids having to take smaller time steps when approaching the minimum voltage.
We fix the time-step size \(\Delta t\) to one minute, for a total of 30 steps, as a representative large step size for which the fully coupled Newton solver remains robust and converges reliably, while still resolving the macroscopic discharge dynamics of interest.

\subsubsection*{Case I: Homogeneous Cubic Cell}

\begin{figure}[htb!]
  \centering
  \begin{tikzpicture}[scale=0.0133,x={(1cm,0cm)}, y={(0.6cm,0.3cm)}, z={(0cm,1cm)}]

  \def\depth{225}       
  \def\width{225}       
  \def\cathode{100}
  \def\separator{25}
  \def\anode{100}
  \def\height{\cathode+\separator+\anode}

  \coordinate (O) at (0,0,0);
  \coordinate (C) at (0,0,\cathode);
  \coordinate (S) at (0,0,\cathode+\separator);
  \coordinate (A) at (0,0,\height);

  \fill[mutedblue!75!gray] (O) -- ++(0,0,\cathode) -- ++(\depth,0,0) -- ++(0,0,-\cathode) -- cycle; 
  \fill[mutedblue!70!gray!95!black] (\depth,0,0) -- ++(0,\width,0) -- ++(0,0,\cathode) -- ++(0,-\width,0) -- cycle; 
  \node at (\depth/2, 0, \cathode/2) {    \begin{tabular}{c}
      negative \\ 
      electrode
    \end{tabular}};

  \fill[mutedgray!20] (C) -- ++(0,0,\separator) -- ++(\depth,0,0) -- ++(0,0,-\separator) -- cycle;
  \fill[mutedgray!25] (\depth,0,\cathode) -- ++(0,\width,0) -- ++(0,0,\separator) -- ++(0,-\width,0) -- cycle;
  \node at (\depth/2, 0, \cathode + \separator/2) {separator};

  \fill[mutedorange!85!white] (S) -- ++(0,0,\anode) -- ++(\depth,0,0) -- ++(0,0,-\anode) -- cycle;
  \fill[mutedorange!90!white] (\depth,0,\cathode+\separator) -- ++(0,\width,0) -- ++(0,0,\anode) -- ++(0,-\width,0) -- cycle;
  \fill[mutedorange!95!white] (A) -- ++(\depth,0,0) -- ++(0,\width,0) -- ++(-\depth,0,0) -- cycle; 
  \node at (\depth/2, 0, \cathode + \separator + \anode/2) {    \begin{tabular}{c}
      positive \\ 
      electrode
    \end{tabular}};

  \draw[<->, thick] (-10,0,0) -- (-10,0,\cathode) node[midway,left] {\SI{100}{\micro\meter}};
  \draw[<->, thick] (-10,0,\cathode) -- (-10,0,\cathode+\separator) node[midway,left] {\SI{25}{\micro\meter}};
  \draw[<->, thick] (-10,0,\cathode+\separator) -- (-10,0,\height) node[midway,left] {\SI{100}{\micro\meter}};

  \draw[<->, thick] (0,0,-10) -- (\depth,0,-10) node[midway,below] {\SI{225}{\micro\meter}};

  \draw[<->, thick] (\depth,0,-10) -- (\depth,\depth,-10) node[midway,below right] {\SI{225}{\micro\meter}};

\end{tikzpicture}
\caption{Diagram of the simulation domain for Cases I and II.}
\label{fig:cubic}
\end{figure}
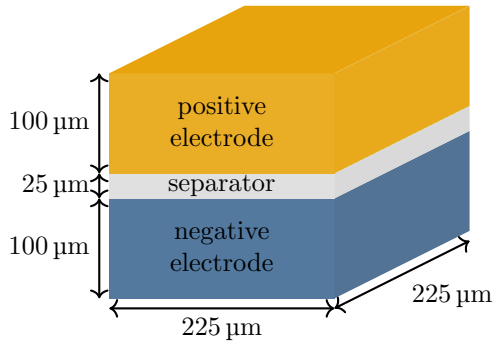

Our first test case is a simple cubic cell with \SI{100}{\micro\meter} thick electrodes and a \SI{25}{\micro\meter} separator as illustrated in \cref{fig:cubic}.
Although not representative of typical operating conditions, the applied current density is prescribed as a two-dimensional Gaussian field in order to induce three-dimensional spatial variation in what would otherwise be a 1D problem.
Let \( I_\mathrm{app} \) be a reference applied current, typically calculated based on the mass loading of the active material.
The spatially varying applied current density 
\( i_{\mathrm{app}} \) on the positive current collector 
\( \Gamma_p \) is given by
\begin{equation}
  i_{\mathrm{app}}(y,z) 
  = \frac{I_{\mathrm{app}}}{\displaystyle 
      \int_{\Gamma_p} 
      \exp\!\left(
        -\frac{(y-y_0)^2}{2\sigma_y^2}
        -\frac{(z-z_0)^2}{2\sigma_z^2}
      \right)  \, \mathrm{d}S }
  \,
  \exp\!\left(
    -\frac{(y-y_0)^2}{2\sigma_y^2}
    -\frac{(z-z_0)^2}{2\sigma_z^2}
  \right),
\end{equation}
where \( (y_0,z_0) \) denotes the center of the collector surface, and
\( \sigma_y \), \( \sigma_z \) are the spreads in the \( y \)-- and \( z \)--directions,
respectively.
The normalization guarantees that
\begin{equation}
  \int_{\Gamma_p} i_{\mathrm{app}}(y,z)  \, \mathrm{d} \mathbf{s}
  \;=\; I_{\mathrm{app}},
\end{equation}
so that the average current density is given by \( I_\mathrm{app} / A\), where \( A = \int_{\Gamma_p}  \, \mathrm{d} S\) is the area of the current collector.
The Gaussian profile is centered in the 
middle of the current collector surface,  
\begin{equation}
  y_0 = \tfrac{1}{2} L_y, 
  \qquad 
  z_0 = \tfrac{1}{2} L_z,
\end{equation}
and the standard deviations are taken as a fixed fraction of the collector
dimensions,  
\begin{equation}
  \sigma_y = 0.1 \, L_y, 
  \qquad 
  \sigma_z = 0.1 \, L_z,
\end{equation}
where \( L_y = L_z = \SI{225}{\micro\meter}\) are the lengths of the collector in the 
\( y \)-- and \( z \)--directions, respectively.

\subsubsection*{Case II: Spatially Varying Properties}

\begin{figure}[htb!]
  \centering
  \begin{subfigure}[t]{0.32\textwidth}
    \centering
    \includegraphics[width=\textwidth]{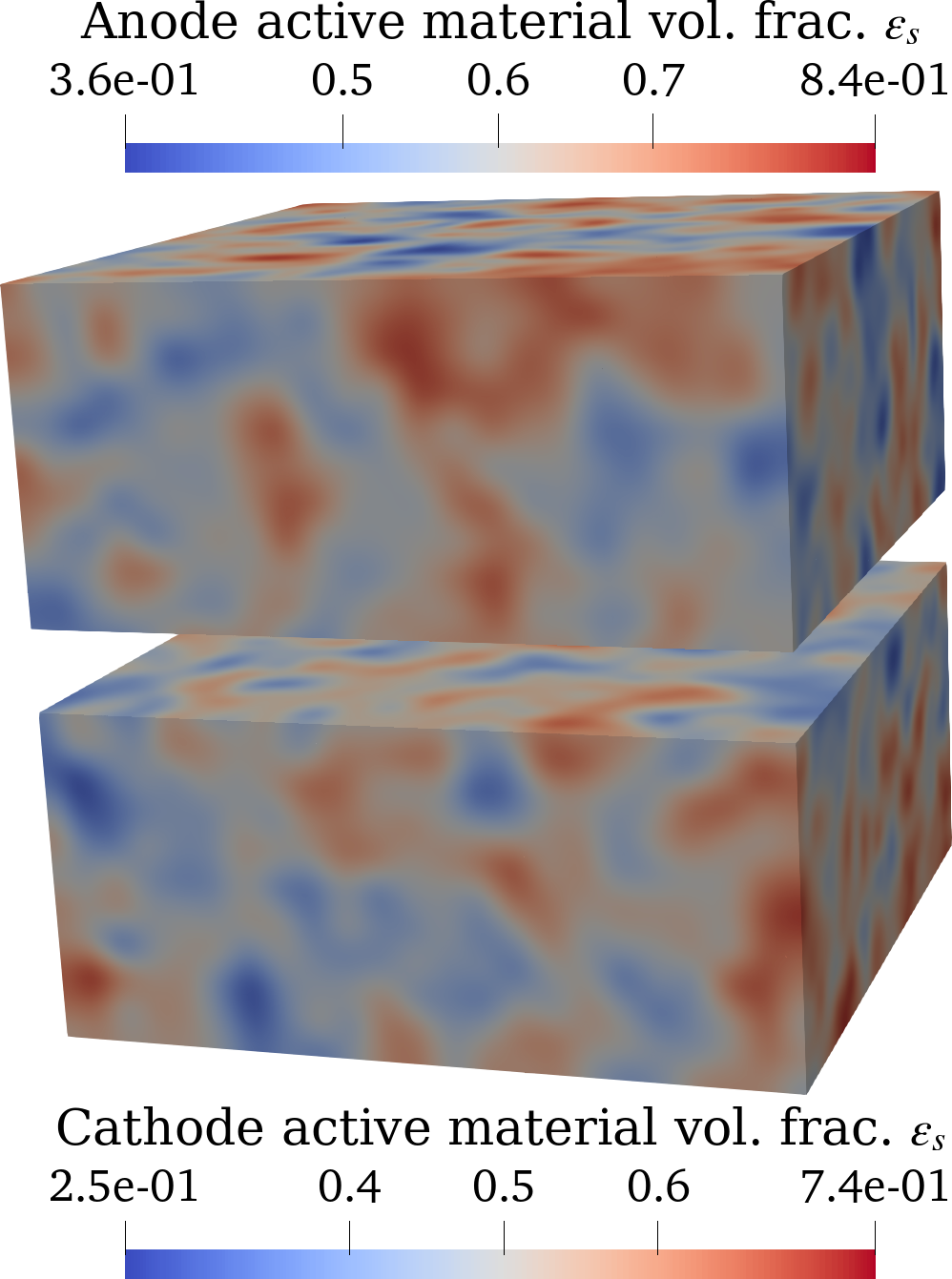}
    \caption{Active material volume fraction.}
    \label{fig:eps_s}
  \end{subfigure}
  \hfill
  \begin{subfigure}[t]{0.32\textwidth}
    \centering
    \includegraphics[width=\textwidth]{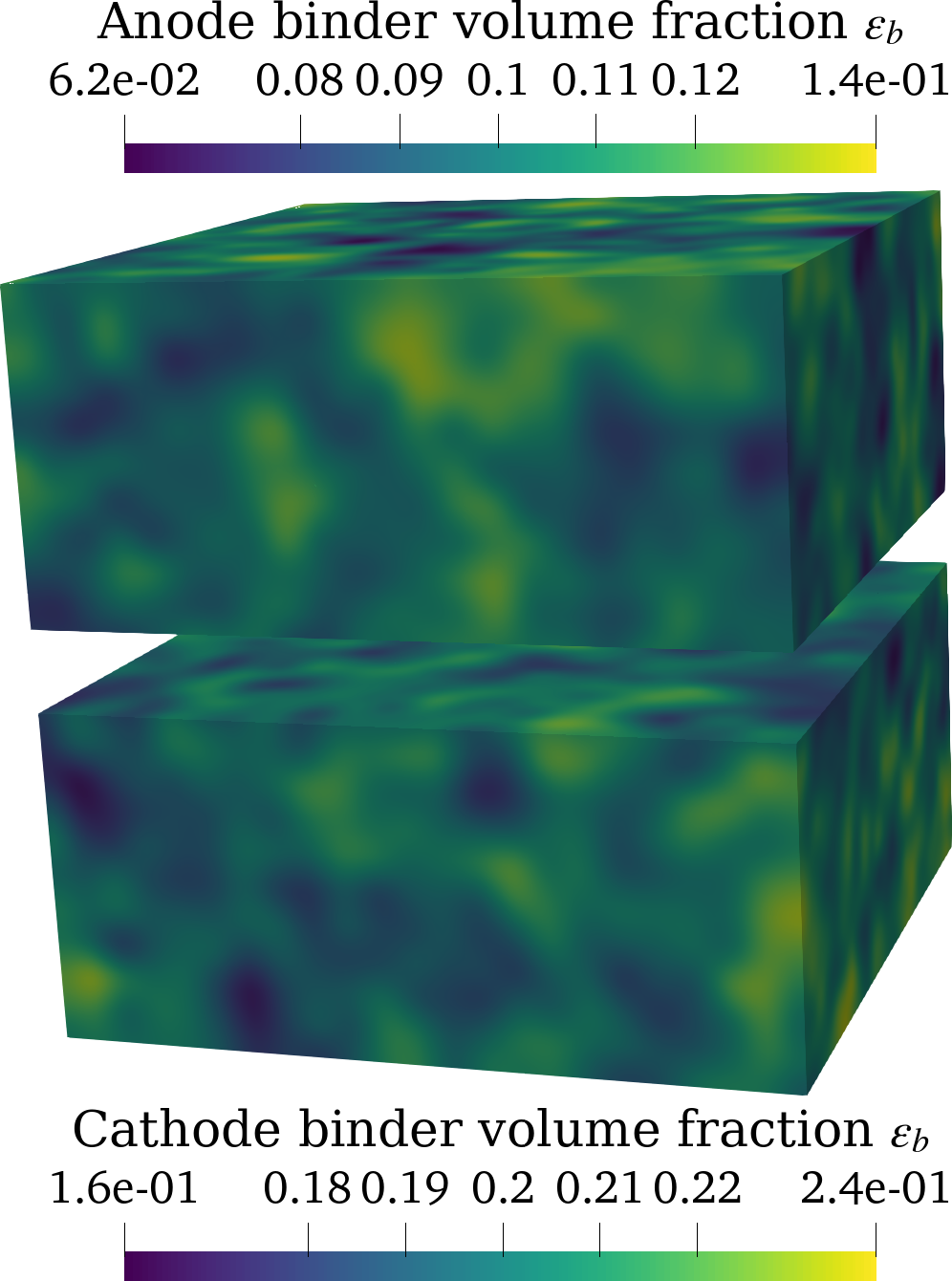}
    \caption{Binder volume fraction.}
    \label{fig:eps_b}
  \end{subfigure}
  \hfill
  \begin{subfigure}[t]{0.32\textwidth}
    \centering
    \includegraphics[width=\textwidth]{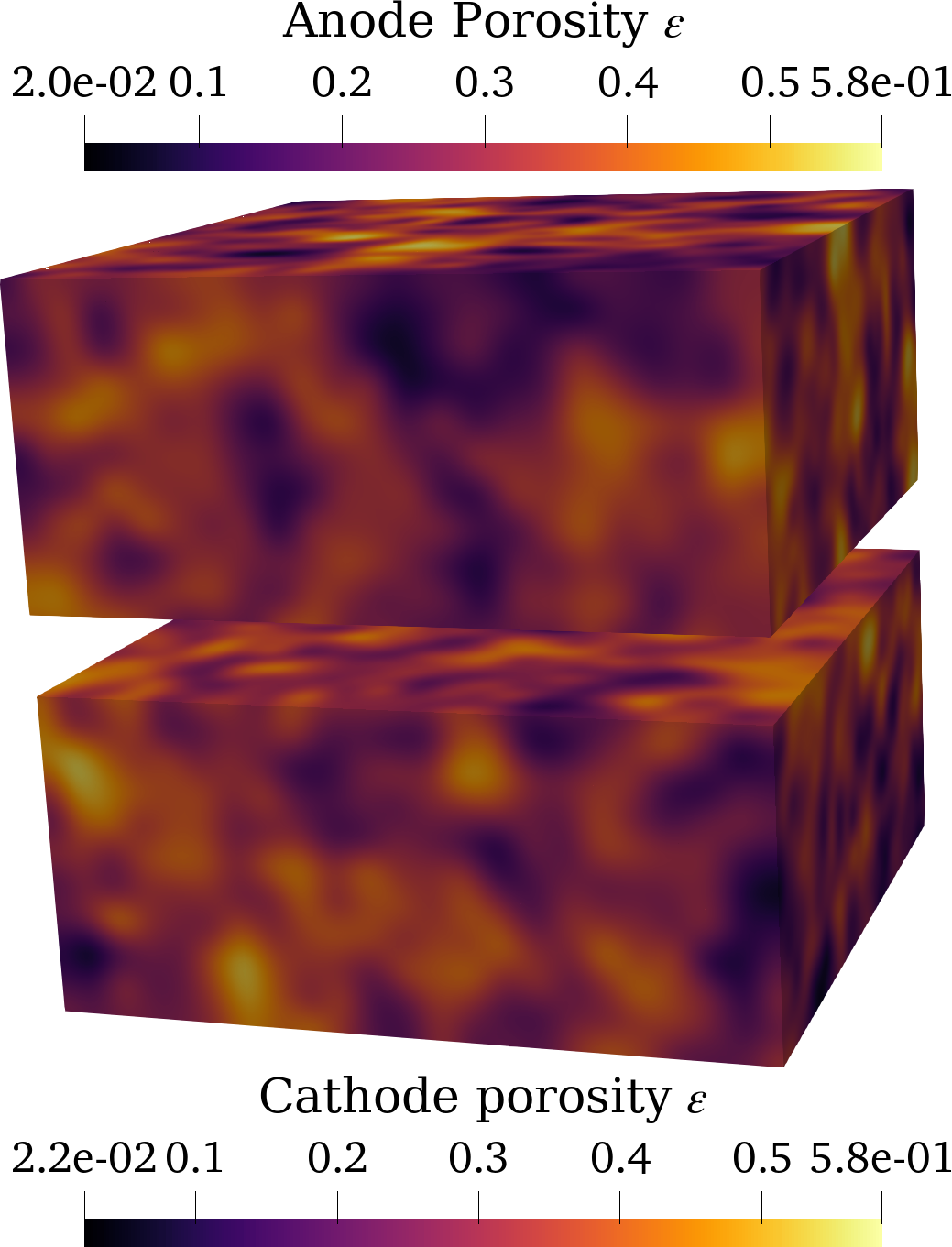}
    \caption{Porosity.}
    \label{fig:porosity}
  \end{subfigure}

  \caption{Spatially varying material properties for Case II.}
  \label{fig:case2}
\end{figure}

Using the same geometry as Case I, we now introduce spatial variations within the domain by adding spatial heterogeneity in material properties.
For the pseudo-2D DFN model, it is typical to consider constant properties since homogenization is already being performed in the other two spatial dimensions.
Typically, spatially-varying properties are captured using simulations where the microstructure (active particles, binder, etc.) is resolved \cite{yan2012three,hutzenlaub2014three,danner2016thick}
However, some electrode-scale heterogeneity (as opposed to particle-scale heterogeneity) can be captured using homogenized porous electrode theory simulations \cite{mistry2020stochasticity, parmananda2022probing}.
In this case, the current is applied uniformly on the current collector face.

To prescribe smooth, spatially correlated heterogeneity for material volume fractions, we construct a filtered random field over the domain $\Omega$.
A coarse, structured grid of resolution $h_c$ is defined on $\Omega$, with nodes $\{(x_i, y_j, z_k)\}_{i,j,k}$.  
On this coarse grid, independent standard Gaussian samples $g_{ijk} \sim \mathcal{N}(0,1)$ are drawn to define a piecewise trilinear interpolant $g(x)$ over $\Omega$.  
To introduce spatial correlation with length scale $h_c$, we apply a reaction-diffusion-type filter~\cite{lazarov2011filters} to obtain a smooth field $\tilde g$ by solving
\[
  (I - h_c^2 \nabla^2)\, \tilde g = g \quad \text{in } \Omega,
\]
on the finer mesh used for the battery simulations, yielding a function $\tilde g (x)$ with reduced high-frequency content.
This is a Mat\'ern type Gaussian random field \cite{croci2018efficient,duswald2024finite}.

\begin{table}[htb!]
\centering
\caption{Prescribed mean and variance values for volume fractions of solid components.}
\label{tab:eps}
\begin{tabular}{lcc}
\hline
\textbf{Component} & \textbf{Mean} ($\bar{\varepsilon}_i$) & \textbf{Variance} ($\sigma_i^2$) \\
\hline
Active material (cathode) & 0.50 & 0.004 \\
Binder (cathode)          & 0.20 & 0.0001 \\
Active material (anode)   & 0.60 & 0.004 \\
Binder (anode)            & 0.10 & 0.0001 \\
\hline
\end{tabular}
\end{table}

The field $\tilde g$ is then normalized to match prescribed statistical moments.  
Let $\bar{\tilde g} = \frac{1}{|\Omega|} \int_\Omega \tilde g (x) \, \mathrm{d} x$ and $\mathrm{Var}[\tilde g] = \frac{1}{|\Omega|} \int_\Omega (\tilde g(x) - \bar{\tilde g})^2 \, \mathrm d x$.  
We define the rescaled field
\[
  \varepsilon_i(x) = \bar{\varepsilon}_i + \sqrt{\frac{\sigma_i^2}{\mathrm{Var}[\tilde g]}} \, (\tilde g (x) - \bar{\tilde g}),
\]
so that $\varepsilon_i(x)$ has mean $\bar{\varepsilon}_i$ and variance $\sigma_i^2$. 
The active-material and binder volume fractions, $\varepsilon_s, \varepsilon_b$, in the cathode and anode regions are generated using the $(\bar{\varepsilon_i}, \sigma_i^2)$ values to given in \cref{tab:eps}.  
The local porosity is then computed as $\varepsilon (x) = 1 - \varepsilon_s(x) - \varepsilon_b(x)$.

\subsubsection*{Case III: Flattened Jelly Roll}
\begin{figure}[htb!]
  \centering
  \begin{subfigure}[t]{0.59\textwidth}
    \centering
    \begin{tikzpicture}[scale=0.02, x={(1cm,0cm)}, y={(0.6cm,0.3cm)}, z={(0cm,1cm)}]

      \def\depth{235}       
      \def\width{130}       
      \def\cathode{65}
      \def\separator{26}
      \def\anode{65}
      \def\cc{13}
      \def\tabw{14}
      \def\tabl{35}
      \def\height{\cathode+\separator+\anode}

      \coordinate (O) at (0,0,0);
      \coordinate (C) at (0,0,\cathode);
      \coordinate (S) at (0,0,\cathode+\separator);
      \coordinate (A) at (0,0,\height);

      \fill[mutedblue!75!gray] (O) -- ++(0,0,\cathode) -- ++(\depth,0,0) -- ++(0,0,-\cathode) -- cycle; 
      \fill[mutedblue!70!gray!95!black] (\depth,0,0) -- ++(0,\width,0) -- ++(0,0,\cathode) -- ++(0,-\width,0) -- cycle; 
      \node at (\depth/2, 0, \cathode/2) {\begin{tabular}{c}
        negative \\ 
        electrode
      \end{tabular}};

      \fill[aluminum!95!black] (0,0,-\cc) -- ++(0,0,\cc) -- ++(\depth,0,0) -- ++(0,0,-\cc) -- cycle; 
      \fill[aluminum!90!black] (\depth,0,-\cc) -- ++(0,\width,0) -- ++(0,0,\cc) -- ++(0,-\width,0) -- cycle; 
      \node[font=\small] at (\depth/2, 0, -\cc/2) {aluminum current collector};


      \fill[black] (\depth-\tabw,0,-\cc) -- ++(\tabw,0,0) -- ++(0,\tabl,0) -- ++(-\tabw,0,0) -- cycle; 
      \fill[white, opacity=0.4] (\depth-\tabw,0,-\cc) -- ++(\tabw,0,0) -- ++(0,\tabl,0) -- ++(-\tabw,0,0) -- cycle; 


      \fill[mutedgray!20] (C) -- ++(0,0,\separator) -- ++(\depth,0,0) -- ++(0,0,-\separator) -- cycle;
      \fill[mutedgray!25] (\depth,0,\cathode) -- ++(0,\width,0) -- ++(0,0,\separator) -- ++(0,-\width,0) -- cycle;
      \node at (\depth/2, 0, \cathode + \separator/2) {separator};

      \fill[mutedorange!85!white] (S) -- ++(0,0,\anode) -- ++(\depth,0,0) -- ++(0,0,-\anode) -- cycle;
      \fill[mutedorange!90!white] (\depth,0,\cathode+\separator) -- ++(0,\width,0) -- ++(0,0,\anode) -- ++(0,-\width,0) -- cycle;
      \fill[mutedorange!95!white] (A) -- ++(\depth,0,0) -- ++(0,\width,0) -- ++(-\depth,0,0) -- cycle; 
      \node at (\depth/2, 0, \cathode + \separator + \anode/2) {\begin{tabular}{c}
        positive \\ 
        electrode
      \end{tabular}};

      \fill[lightcopper!95!black] (A) -- ++(0,0,\cc) -- ++(\depth,0,0) -- ++(0,0,-\cc) -- cycle;
      \fill[lightcopper!90!black] (\depth,0,\cathode+\separator+\anode) -- ++(0,\width,0) -- ++(0,0,\cc) -- ++(0,-\width,0) -- cycle;
      \fill[lightcopper!85!black] (0,0,\height+\cc) -- ++(\depth,0,0) -- ++(0,\width,0) -- ++(-\depth,0,0) -- cycle; 
      \node[font=\small] at (\depth/2, 0, \cathode + \separator + \anode + \cc/2) {copper current collector};

      \fill[black] (0,\width-\tabl,\height+\cc) -- ++(\tabw,0,0) -- ++(0,\tabl,0) -- ++(-\tabw,0,0) -- cycle; 
      \node[font=\footnotesize] at (\tabw+20.,\width-\tabl/2,\height+\cc) {$\leftarrow$tab};

      \draw[<->, thick] (-10,0,0) -- (-10,0,\cathode) node[midway,left] {\SI{50}{\micro\meter}};
      \draw[<->, thick] (-10,0,\cathode) -- (-10,0,\cathode+\separator) node[midway,left] {\SI{20}{\micro\meter}};
      \draw[<->, thick] (-10,0,\cathode+\separator) -- (-10,0,\height) node[midway,left] {\SI{50}{\micro\meter}};
      \draw[<->, thick] (-10,0,\height) -- (-10,0,\height+\cc) node[midway,left] {\SI{10}{\micro\meter}};
      \draw[<->, thick] (-10,0,-\cc) -- (-10,0,0) node[midway,left] {\SI{10}{\micro\meter}};

      \draw[<->, thick] (0,0,-10-\cc) -- (\depth,0,-10-\cc) node[midway,below] {\SI{1}{\meter}};

      \draw[<->, thick]
        (\depth,0,-10-\cc) -- node[midway, sloped, below]{\SI{5}{\centi\meter}}
        node[pos=0.55, sloped, above=2pt, font=\footnotesize]{$\leftarrow$tab}
        (\depth,\width,-10-\cc);

    \end{tikzpicture}
    \caption{Diagram of a flattened jelly roll.}
    \label{fig:flattened}
  \end{subfigure}
  \hfill
  \begin{subfigure}[t]{0.39\textwidth}
    \centering
    \begin{tikzpicture}[scale=0.80]

    \def\rzero{0.2}        
    \def\pitch{0.12}       
    \def\thick{0.22}       
    \def\turns{6.5*pi}     
    \def\samples{500}
    \def\tabwidth{0.05*pi} 

    \colorlet{anode}{mutedorange!90!white}
    \colorlet{cathode}{mutedblue!80!gray}
    \colorlet{separator}{mutedgray!20}
    \colorlet{al}{aluminum!90!black}
    \colorlet{cu}{lightcopper!90!black}
    \colorlet{tab}{black}

    \newcommand{\spiralband}[5]{%
      \path[fill=#1, draw=none]
        plot[domain=#4:#5, samples=\samples, variable=\t, smooth]
          ({(\rzero+\pitch*\t+#2)*cos(\t r)}, {(\rzero+\pitch*\t+#2)*sin(\t r)})
        --
        plot[domain=#5:#4, samples=\samples, variable=\t, smooth]
          ({(\rzero+\pitch*\t+#3)*cos(\t r)}, {(\rzero+\pitch*\t+#3)*sin(\t r)})
        -- cycle;
    }

    \spiralband{al}{0.0*\thick}{0.3*\thick}{0}{\turns}        
    \spiralband{cathode}{0.3*\thick}{1.5*\thick}{0}{\turns}   
    \spiralband{separator}{1.5*\thick}{1.9*\thick}{0}{\turns} 
    \spiralband{anode}{1.9*\thick}{3.1*\thick}{0}{\turns}     
    \spiralband{cu}{3.1*\thick}{3.4*\thick}{0}{\turns}        


    \spiralband{tab}{-0.4*\thick}{0.0*\thick}{0}{10.0*\tabwidth}

    \spiralband{tab}{3.4*\thick}{3.7*\thick}{\turns-\tabwidth}{\turns}

  \node[font=\footnotesize, anchor=east, text=al!70!black] at (-1.1,2.68) {Al CC};
  \draw[-{Stealth[length=2.5mm,width=1.0mm]}, color=al!70!black, thick]
    (-1.1,2.68) -- (0,2.68);

  \node[font=\footnotesize, anchor=east, text=cathode!70!black] at (-1.1,2.975) {cathode};
  \draw[-{Stealth[length=2.5mm,width=1.0mm]}, color=cathode!70!black, thick]
    (-1.1,2.975) -- (0,2.85);

  \node[font=\footnotesize, anchor=east, text=separator!70!black] at (-1.1,3.20) {separator};
  \draw[-{Stealth[length=2.5mm,width=1.0mm]}, color=separator!70!black, thick]
    (-1.1,3.20) -- (0,3.02);

  \node[font=\footnotesize, anchor=east, text=anode!70!black] at (-1.1,3.455) {anode};
  \draw[-{Stealth[length=2.5mm,width=1.0mm]}, color=anode!70!black, thick]
    (-1.1,3.455) -- (0,3.2);

  \node[font=\footnotesize, anchor=east, text=cu!70!black] at (-1.1,3.69) {Cu CC};
  \draw[-{Stealth[length=2.5mm,width=1.0mm]}, color=cu!70!black, thick]
    (-1.1,3.69) -- (0,3.36);

    \node[font=\footnotesize, color=black]     at (-0.2,0.1)   {tab};
    \node[font=\footnotesize, color=black]     at (0.4,3.6) {tab};

  \end{tikzpicture}
    \caption{Diagram of a jelly roll cross-section.}
    \label{fig:jellycross}
  \end{subfigure}
\caption{Diagrams of the simulation domain of a flattened jelly roll for Case III. The diagrams are not to scale (note the dimensions in (a)). The curvature in (b) is for illustrative purposes only and is not considered in the simulations.}
\label{fig:jelly}
\end{figure}
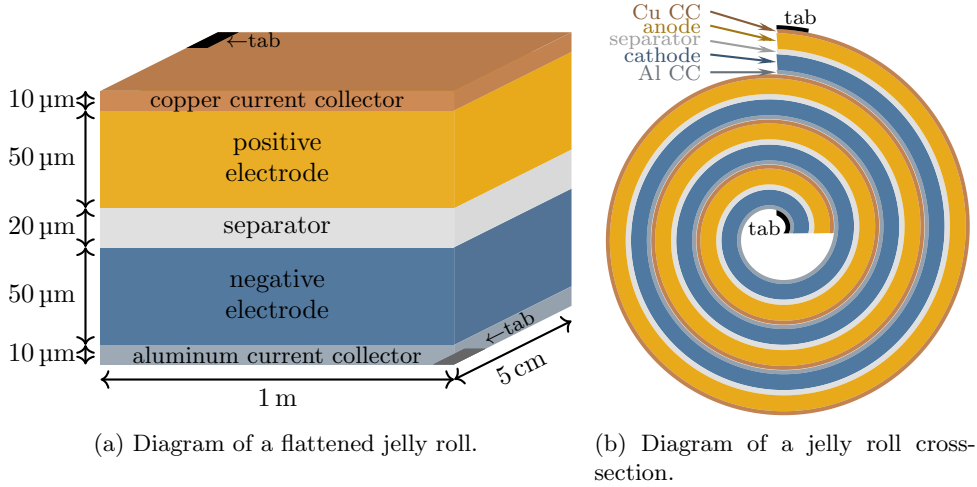

A typical cylindrical lithium-ion cell employs a \emph{jelly-roll} architecture, formed by winding stacked electrode, separator, and current collector layers around a central mandrel, as illustrated in \cref{fig:jellycross}.
In this work, we use a flattened jelly-roll geometry obtained by unwrapping the spiral into a planar stack (\cref{fig:flattened}).
Since the stack thickness is negligible compared to the radius of curvature, neglecting curvature has little impact on the relevant electrochemical behavior.
For simplicity, the model features single-sided electrodes rather than the more compact double-sided configuration commonly used in practical cells.

The dimensions of the simulation domain are shown in \cref{fig:flattened}, not to scale.
The simulation domain is quite anisotropic: \SI{1e-4}{\meter} $\times$ \SI{0.05}{\meter} $\times$ \SI{1}{\meter}.
The resulting hexahedral mesh used for strong scaling studies is also very anisotropic as mentioned in \cref{sec:strong}.

For this case, we explicitly model the current collectors, rather than treating them as boundary surfaces as in the other cases, since the applied current is injected only through localized tabs.
In these subdomains, only the solid-phase potential is solved \cref{eq:phis}, with zero porosity and high solid-phase conductivity, as noted in \cref{sec:model_parameters}.
The applied current and grounded potential boundary conditions \cref{eq:Neumann,eq:Dirichlet} are imposed on tabs of dimensions \SI{5}{\milli\meter} $\times$ \SI{10}{\milli\meter} located on the current collectors.
These tabs are positioned on opposite ends of the domain, corresponding to the top/bottom and interior/exterior of the rolled cell.
Despite the high conductivity of the current collectors, their finite dimensions lead to nonuniform current distribution and consequently three-dimensional spatial variations.

Severe anisotropy can affect the performance of AMG.
Therefore, for this case, more conservative parameters are used for BoomerAMG (e.g. less aggressive coarsening), as detail in \cref{sec:solver_parameters}.
In addition, we increase the GMRES restart size from the default 30 to 100.

\subsubsection*{Case IV: Interpenetrating Gyroid}

\begin{figure}
  \centering
  \includegraphics[width=0.5\linewidth]{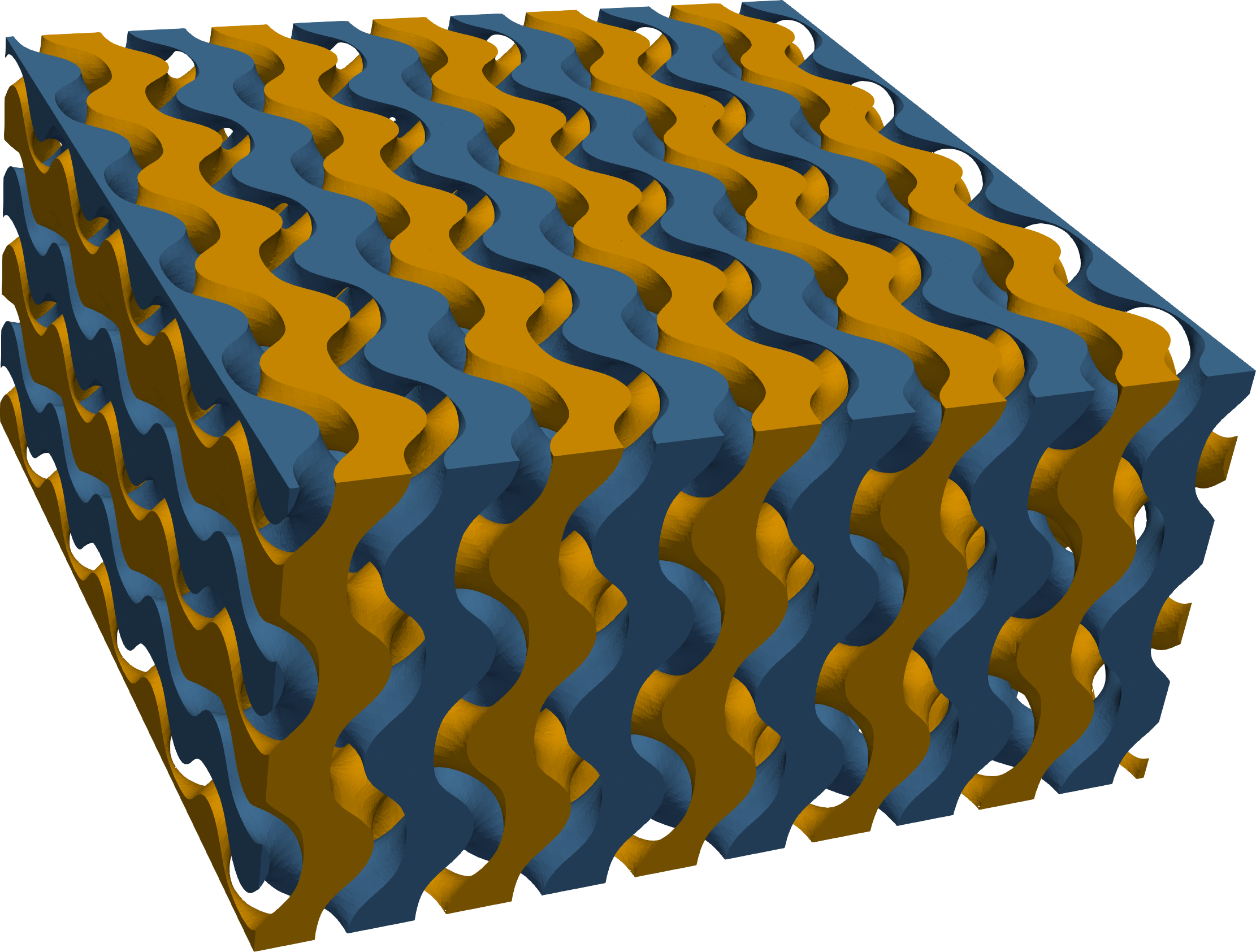}
  \caption{Simulation domain for Case IV with TPMS gyroid electrodes where the negative electrode (anode) is orange, the positive electrode (cathode) is blue and the void in between is filled with electrolyte.}
  \label{fig:tpms}
\end{figure}

Interpenetrating electrode architectures have recently attracted increasing attention in the energy storage literature due to their high surface-area-to-volume ratios and tunable transport pathways, leading to thick cells that are not limited by ion transport \cite{xue2024interpenetrated,xue2026architecting,wiesner2023additive,wang2026ultra}.

We generate triply-periodic minimal surface (TPMS) electrode geometries in a physical computational domain
\(\Omega = [0,L_x]\times[0,L_y]\times[0,L_z]\), with
\((L_x,L_y,L_z) = (\SI{1}{mm}, \SI{1}{mm}, \SI{0.5}{mm})\),
by constructing signed distance fields derived from the gyroid surface representation.
Introducing dimensionless coordinates
\(
\hat x = x/L_x,\;
\hat y = y/L_y,\;
\hat z = z/L_z,
\)
the phase-shifted gyroid level-set functions for the two electrodes are
\begin{align}
  s_n(\hat x,\hat y,\hat z) &=
  \sin(\pi \hat x)\cos(\pi \hat y)
  + \sin(\pi \hat y)\cos(\pi \hat z)
  + \sin(\pi \hat z)\cos(\pi \hat x) - 1.3, \\
  s_p(\hat x,\hat y,\hat z) &=
  \sin(\pi \hat x)\cos(\pi \hat y + \pi)
  + \sin(\pi \hat y + \pi)\cos(-\pi \hat z)
  + \sin(-\pi \hat z)\cos(\pi \hat x) - 1.3.
\end{align}
Corresponding indicator regions are defined from the signed distance inequalities
\begin{equation}
  |s_n| \le \frac{t}{2}, \qquad
  |s_p| \le \frac{t}{2},
\end{equation}
with \(t = 1.104\), yielding disjoint domains representing the anode (\(\Omega_n\)) and cathode (\(\Omega_p\)) phases.

To obtain a mesh conforming to the electrode--electrolyte interfaces, we solve a Target-Matrix Optimization Paradigm (TMOP) problem~\cite{dobrev2019target} using the MFEM framework~\cite{anderson2021mfem}.
This optimization procedure, frequently applied in level-set topology optimization, minimizes mesh distortion while aligning element faces with material interfaces~\cite{schmidt2024level}.
This is a larger version of the electrodes and mesh used in \cite{cross2024viability}.

\subsection{Impact of system ordering on block Gauss-Seidel preconditioner}
\label{sec:ordering}
As noted in \cref{sec:preconditioning}, there are 24 different orderings of the $4\times 4$ block system and thus 24 different choices for the block Gauss-Seidel preconditioner, $P_{GS}$ in \cref{eq:PGS}.
We tested all $24$ orderings for the four cases and measured the Krylov iteration counts for the 200 processor simulations used in the strong scaling study in \cref{sec:strong}.

Across all tests, the difference between the best and worst orderings was modest: iteration counts varied by only $6$--$8\%$, and the top five or six orderings differed by less than $2\%$.
Two permutations, $(\phi_e,\, c_s,\, \phi_s,\, c_e)$ and $(c_e,\, \phi_e,\, c_s,\, \phi_s)$, consistently ranked first or second across all geometries, indicating that the most effective orderings are robust with respect to mesh structure and local heterogeneity.
All weak- and strong-scaling results reported below therefore use the $(\phi_e,\, c_s,\, \phi_s,\, c_e)$ ordering.
Because even the least favorable ordering increases the iteration count by only about $6$--$8\%$, the overall performance of the multiplicative block preconditioner is not highly sensitive to the ordering choice.
Note that the optimal ordering could differ if material parameters or operating conditions substantially modify the relative strength of the cross-couplings, which may affect the performance of some orderings.

\subsection{Weak scaling}
\label{sec:weak}

\begin{figure}[htb!]
\centering
\begin{tikzpicture}
    \begin{semilogxaxis}[
        log ticks with fixed point,
        scaled ticks=false,
        tick scale binop=\times,
        width=3.75in,
        xlabel={number of processors},
        ylabel={solve time (min)},
        title={Weak Scaling (Cubic Cell, $\sim 142$k DoFs/processor)},
        xtick={25,50,100,200,400,800,1600},
        xticklabels={25,50,100,200,400,800,1600},
        log basis x=2,
        grid=both,
        grid style={dotted, gray!30},
        legend style={at={(0.98,0.02)}, anchor=south east, font=\small},
        ]

        \addplot[thick,color=blue,mark=*] coordinates {
            (25,  2.132)
            (50,  2.703)
            (100, 3.964)
            (200, 4.239)
            (400, 4.468)
            (800, 4.602)
            (1600,4.969)
        };

        \addplot[thick,color=blue,mark=*,dashed] coordinates {
            (25,  0.748)
            (50,  1.094)
            (100, 1.869)
            (200, 2.019)
            (400, 2.226)
            (800, 2.261)
            (1600,2.458)
        };

        \node[blue, font=\small] at (axis cs:25,  1.9)  {(1620)};
        \node[blue, font=\small] at (axis cs:52,  2.35)  {(1706)};
        \node[blue, font=\small] at (axis cs:115, 3.65)  {(1760)};
        \node[blue, font=\small] at (axis cs:200, 3.90)  {(1818)};
        \node[blue, font=\small] at (axis cs:400, 4.1)  {(1944)};
        \node[blue, font=\small] at (axis cs:800, 4.22)  {(1904)};
        \node[blue, font=\small] at (axis cs:1600,4.6)  {(1968)};

        \addplot[thick,color=orange,mark=diamond*] coordinates {
            (25,  2.461)
            (50,  3.191)
            (100, 4.314)
            (200, 4.619)
            (400, 4.965)
            (800, 5.057)
            (1600,5.436)
        };

        \addplot[thick,color=orange,mark=diamond*,dashed] coordinates {
            (25,  0.712)
            (50,  1.086)
            (100, 1.921)
            (200, 2.106)
            (400, 2.393)
            (800, 2.411)
            (1600,2.628)
        };

        \node[orange, font=\small] at (axis cs:22,  2.8)  {(1180)};
        \node[orange, font=\small] at (axis cs:40,  3.5)  {(1287)};
        \node[orange, font=\small] at (axis cs:100,4.60)  {(1338)};
        \node[orange, font=\small] at (axis cs:200,5.0)  {(1419)};
        \node[orange, font=\small] at (axis cs:400,5.25)  {(1568)};
        \node[orange, font=\small] at (axis cs:790,5.35)  {(1542)};
        \node[orange, font=\small] at (axis cs:1600,5.6) {(1614)};

        \legend{BJ, BJ (linear solve), BGS, BGS (linear solve)};
    \end{semilogxaxis}
\end{tikzpicture}
\caption{Weak scaling performance with 30 time steps and approximately 142k DoFs per processor for Case~I. The total number of GMRES iterations is given in parentheses.}
\label{fig:weak_scaling}
\end{figure}
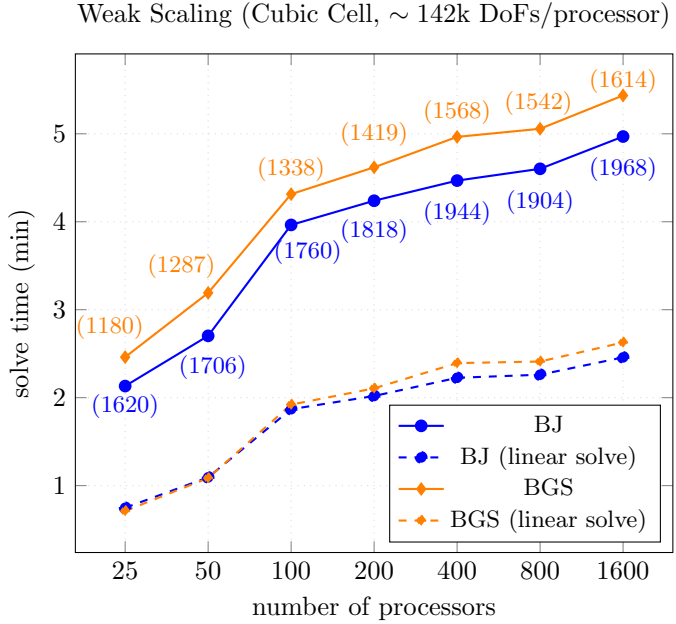

We investigate the weak scalability of the BJ and BGS preconditioners by keeping the local problem size fixed at approximately $1.4\times 10^5$ DoFs per processor while increasing the processor count.
Since mesh refinement is simpler for Case I, we use this case for the weak scaling study.
Figure~\ref{fig:weak_scaling} reports the total nonlinear solve time, the time spent in the linear solve, and the corresponding GMRES iteration counts.
For both BJ and BGS, the GMRES iteration counts increase only mildly between 25 and 1600 processors (approximately $36\%$ for BJ and $21\%$ for BGS), indicating that the preconditioners retain their effectiveness as the global problem size grows.
The remaining increases in runtime are primarily due to communication overhead.
The jump between 50 and 100 processors corresponds to moving from one to two sockets on the same node.
The increase from 100 to 200 processors reflects the transition from intra-node to inter-node communication.
Beyond this point, the total solve time grows only modestly from 200 to 1600 processors (approximately $17\%$), demonstrating near-ideal weak scaling for both preconditioners in this regime.
The BJ preconditioner is consistently faster in total time, although BGS achieves slightly lower linear iteration counts.

\subsection{Strong scaling}
\label{sec:strong}

For each test case, we now perform a strong scaling study, i.e. the problem size is fixed and the number of processors is increased.
The number of processors starts at 200 (2 nodes), and is doubled until 1600 processors (16 nodes).
The problem size is chosen to be close to the maximum size allowed by memory requirements for the 200 processor case.
Note that, at this scale, direct solvers are unfeasible due to memory requirements, so we can only test iterative solvers.

\begin{table}[htb!]
\centering
\caption{Summary of meshes and degrees of freedom for all test cases.}
\label{tab:meshes}
\begin{tabular}{c c S[table-format=1.1e2] S[table-format=1.1e2] c S[table-format=1.2e2] S[table-format=1.2e2]}
\hline
Case & Mesh & {Elements} & {Electrode DoFs} & $N_c$ & {Particle DoFs} & {Total DoFs} \\
\hline
I   & Hex & 8.7e6 & 2.6e7 & 10 & 8.7e7 & 1.13e8 \\
II  & Hex & 8.7e6 & 2.6e7 & 10 & 8.7e7 & 1.13e8 \\
III & Hex & 8.5e6 & 2.6e7 & 10 & 8.5e7 & 1.12e8 \\
IV  & Tet & 3.6e7 & 1.8e7 & 6  & 2.18e8 & 2.37e8 \\
\hline
\end{tabular}
\end{table}

Table~\ref{tab:meshes} summarizes the meshes and degrees of freedom for all test cases, including the particle-level discretization used in each case.
For Cases I-III, hexahedral elements are used, with Cases I and II using identical meshes.
The elements are equi-spaced in each direction.
For Cases I and II, there are 207 elements in the direction of the cell thickness, so that elements are split perfectly across the positive electrode, separator, and negative domains.
The other directions are 205 element-wide, leading to just under 8.7 million elements, so over 26 million electrode-level degrees of freedom.
In the particle-direction, we use $N_c=10$ points.
Since in our implementation these are element-wise unknowns, we get just under 87 million particle-level degrees of freedom, for a total of over 113 million degrees of freedom.
For Case III, we use 56 elements in the direction of cell thickness, 85 elements for the cell height, and 1800 elements along the length of the roll.
This leads to over 8.5 million elements, so over 26 million electrode-level degrees of freedom.
Again, we pick $N_c=10$, leading to 85 million particle-level degrees of freedom, and thus 112 million total degrees of freedom.

For Case IV, a tetrahedral mesh with over 36 million elements is used, leading to over 18 million electrode-level degrees of freedom.
For this case, we choose a smaller $N_c=6$ number of points in the particle direction.
However, since these are element-wise unknowns, we get around 218 million particle-level degrees of freedom, for a total of 237 million degrees of freedom.
Recall from \cref{sec:preconditioning}, that the particle system is tridiagonal, so not contributing that many non-zeros compared to the electrode-level system.

    
    

\begin{figure}[htb!]
\centering
\begin{tikzpicture}
    \begin{loglogaxis}[
      log ticks with fixed point,
      scaled ticks=false,
      tick scale binop=\times,
        width=3.75in,
        xlabel={number of processors},
        ylabel={solve time (min)},
        title={Strong Scaling (Cubic Cell)},
        xtick={200,400,800,1600},
        ytick={20,10,5,2.5,1.25},
        log basis x=2,
        log basis y=2,
        grid=both,
        grid style={dotted, gray!30},
        legend style={at={(0.98,0.98)}, anchor=north east},
        ]

        \addplot[thick,color=blue,mark=*] coordinates {
            (200, 17.13)
            (400, 9.31)
            (800, 4.68)
            (1600, 2.77)
            }
            node[pos=0.0, below, yshift=-1mm]{(1870)}
            node[pos=0.33, below, yshift=-2mm]{(2121)}
            node[pos=0.67, below, yshift=-2mm]{(1925)}
            node[pos=1.0, below, yshift=-1.2mm]{(1898)}
            ;

        \addplot[thick,color=blue,mark=*,dashed] coordinates {
            (200, 9.59)
            (400, 5.18)
            (800, 2.27)
            (1600, 1.25)
            };

        \addplot[thick,color=orange,mark=diamond*] coordinates {
            (200, 19.89)
            (400, 10.58)
            (800, 5.30)
            (1600, 3.06)
            }
            node[pos=0.0, above]{(1514)}
            node[pos=0.33, above, xshift=0.9mm]{(1816)}
            node[pos=0.67, above, xshift=0.9mm]{(1581)}
            node[pos=1.0, above, xshift=0.4mm]{(1546)}
            ;

        \addplot[thick,color=orange,mark=diamond*,dashed] coordinates {
            (200, 11.06)
            (400, 5.81)
            (800, 2.55)
            (1600, 1.37)
            };

        \addplot[thick,color=black,dotted,opacity=0.5] coordinates {
            (200, 17.13)
            (400, 8.56)
            (800, 4.28)
            (1600, 2.14)
        };
        \addplot[thick,color=black,dotted,opacity=0.5] coordinates {
            (200, 9.59)
            (400, 4.79)
            (800, 2.40)
            (1600, 1.20)
        };
        \addplot[thick,color=black,dotted,opacity=0.5] coordinates {
            (200, 19.89)
            (400, 9.95)
            (800, 4.97)
            (1600, 2.49)
        };
        \addplot[thick,color=black,dotted,opacity=0.5] coordinates {
            (200, 11.06)
            (400, 5.53)
            (800, 2.77)
            (1600, 1.38)
        };

        \legend{BJ, BJ (linear solve), BGS, BGS (linear solve)};
    \end{loglogaxis}
    \end{tikzpicture}
    \caption{Strong scaling timings for Case I, the homogeneous cubic cell problem (113M DoFs). The total number of GMRES iterations are given in parentheses. Ideal scaling is given by dotted lines.}
    \label{fig:scaling}
\end{figure}
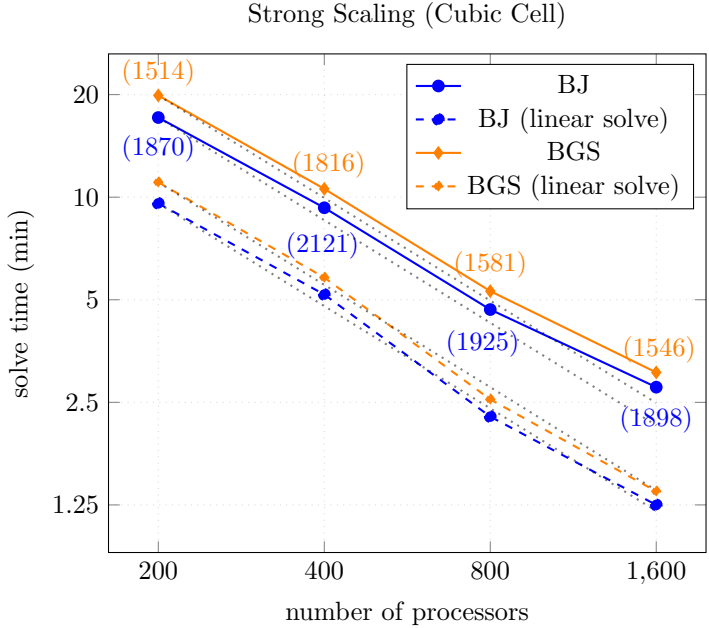

For Case~I, the strong scaling results are shown in \cref{fig:scaling}.
Both BJ and BGS exhibit substantial reductions in solve time as the processor count increases from 200 to 1600.
For BJ, the total solve time decreases from approximately \(1028\) seconds at 200 processors to \(166\) seconds at 1600 processors, corresponding to a \(6.2\times\) speedup and a parallel efficiency of \(0.78\) relative to ideal \(8\times\) scaling.
BGS shows a similar trend, decreasing from \(1193\) seconds to \(184\) seconds over the same range, yielding a \(6.5\times\) speedup and an efficiency of \(0.81\).
The linear solve times follow the same pattern, and the GMRES iteration counts remain nearly constant (variations below \(15\%\) for BJ and \(20\%\) for BGS), indicating that the preconditioner performance is relatively robust under mesh partitioning.
BJ consistently attains lower wall-clock times despite slightly higher iteration counts, while BGS incurs a higher per-iteration cost.
The total solve time increases at a faster rate from 800 to 1600 processors, indicating higher communication costs.
Note that at 800 processors, there are around 142k DoFs per processor, which was found to be a good balance in the weak scaling study.
At 71k DoFs per processor, communication is starting to become a bottleneck.

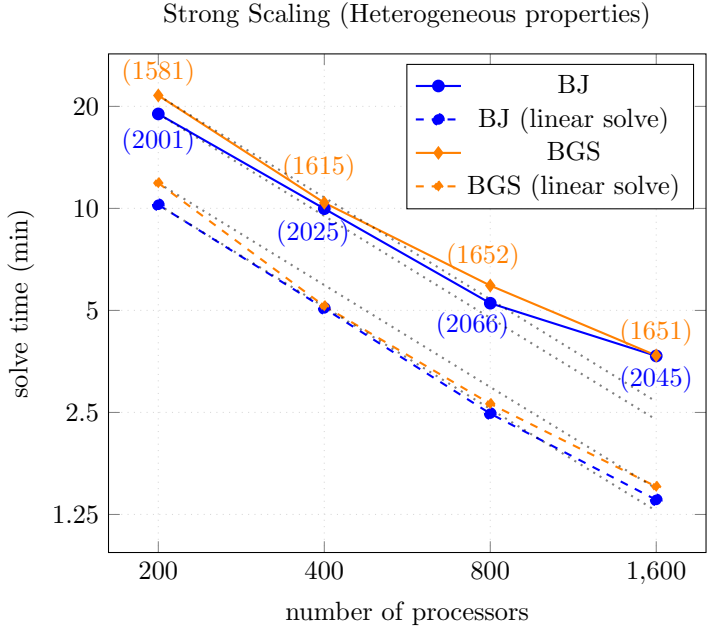
\begin{figure}[htb!]
\centering
\begin{tikzpicture}
    \begin{loglogaxis}[
      log ticks with fixed point,
      scaled ticks=false,
      tick scale binop=\times,
        width=3.75in,
        xlabel={number of processors},
        ylabel={solve time (min)},
        title={Strong Scaling (Heterogeneous properties)},
        xtick={200,400,800,1600},
        ytick={20,10,5,2.5,1.25},
        log basis x=2,
        log basis y=2,
        grid=both,
        grid style={dotted, gray!30},
        legend style={at={(0.98,0.98)}, anchor=north east},
        ]

        \addplot[thick,color=blue,mark=*] coordinates {
            (200, 19.00)
            (400, 9.98)
            (800, 5.25)
            (1600, 3.67)
            }
            node[pos=0.0, below, yshift=-0.6mm]{(2001)}
            node[pos=0.33, below, yshift=-1.1mm]{(2025)}
            node[pos=0.67, below, yshift=-1.2mm]{(2066)}
            node[pos=1.0, below]{(2045)}
            ;

        \addplot[thick,color=blue,mark=*,dashed] coordinates {
            (200, 10.25)
            (400, 5.07)
            (800, 2.48)
            (1600, 1.38)
            };

        \addplot[thick,color=orange,mark=diamond*] coordinates {
            (200, 21.54)
            (400, 10.39)
            (800, 5.92)
            (1600, 3.68)
            }
            node[pos=0.0, above]{(1581)}
            node[pos=0.33, above, xshift=1.5mm, yshift=0.2mm]{(1615)}
            node[pos=0.67, above, xshift=0.6mm]{(1652)}
            node[pos=1.0, above]{(1651)}
            ;

        \addplot[thick,color=orange,mark=diamond*,dashed] coordinates {
            (200, 11.88)
            (400, 5.15)
            (800, 2.65)
            (1600, 1.51)
            };

        \addplot[thick,color=black,dotted,opacity=0.5] coordinates {
            (200, 19.00)
            (400, 9.50)
            (800, 4.75)
            (1600, 2.38)
        };
        \addplot[thick,color=black,dotted,opacity=0.5] coordinates {
            (200, 10.25)
            (400, 5.12)
            (800, 2.56)
            (1600, 1.28)
        };
        \addplot[thick,color=black,dotted,opacity=0.5] coordinates {
            (200, 21.54)
            (400, 10.77)
            (800, 5.39)
            (1600, 2.69)
        };
        \addplot[thick,color=black,dotted,opacity=0.5] coordinates {
            (200, 11.88)
            (400, 5.94)
            (800, 2.97)
            (1600, 1.49)
        };

        \legend{BJ, BJ (linear solve), BGS, BGS (linear solve)};
    \end{loglogaxis}
    \end{tikzpicture}
    \caption{Strong scaling timings for Case II, the cubic cell with heterogeneous properties (113M DoFs). The total number of GMRES iterations are given in parentheses. Ideal scaling is given by dotted lines.}
    \label{fig:scaling_gaussian}
\end{figure}

For Case~II, the strong scaling results for the heterogeneous electrodes are shown in \cref{fig:scaling_gaussian}.  
Compared with the homogeneous Case~I, both BJ and BGS typically require on the order of \(5\%\) more GMRES iterations, reflecting a mild increased difficulty introduced by spatial variations in the active-material and binder volume fractions.  
Nevertheless, both preconditioners continue to exhibit strong parallel scalability.
For the BJ preconditioner, the total solve time decreases from approximately \(1140\)~seconds at 200 processors to \(220\)~seconds at 1600 processors, corresponding to a \(5.2\times\) speedup and a parallel efficiency of \(0.65\).  
The BGS preconditioner shows a similar reduction, decreasing from \(1292\)~seconds to \(221\)~seconds over the same range, yielding a \(5.9\times\) speedup and an efficiency of \(0.73\).  
The linear solve times exhibit nearly identical trends, and iteration counts remain essentially flat across processor counts, confirming that the preconditioners remain robust despite the introduction of spatially correlated heterogeneity.  
As in Case~I, BJ achieves the lowest wall-clock times, but BGS maintains comparable parallel efficiency and slightly lower iteration counts.
Again, communication starts to become a bottleneck at 1600 processors.

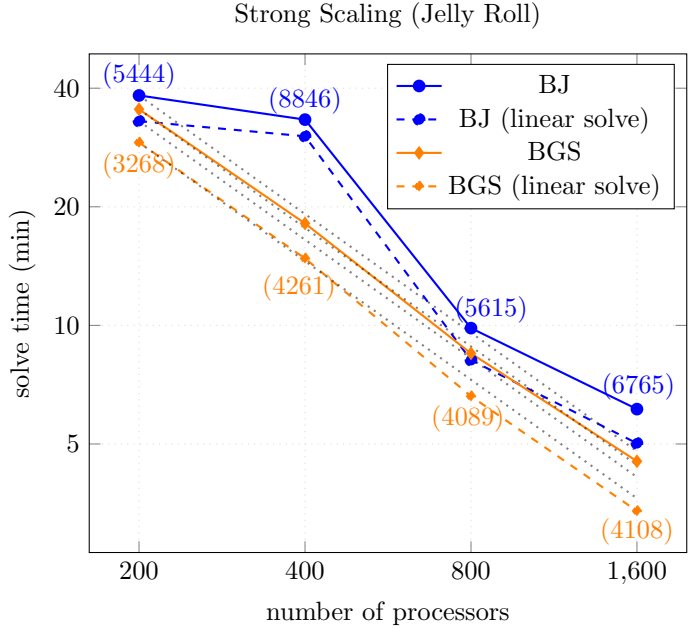
\begin{figure}[htb!]
\centering
\begin{tikzpicture}
    \begin{loglogaxis}[
      log ticks with fixed point,
      scaled ticks=false,
      tick scale binop=\times,
        width=3.75in,
        xlabel={number of processors},
        ylabel={solve time (min)},
        title={Strong Scaling (Jelly Roll)},
        xtick={200,400,800,1600},
        ytick={5,10,20,40},
        log basis x=2,
        log basis y=2,
        grid=both,
        grid style={dotted, gray!30},
        legend style={at={(0.98,0.98)}, anchor=north east},
        ]

        \addplot[thick,color=blue,mark=*] coordinates {
            (200, 38.34)
            (400, 33.29)
            (800, 9.84)
            (1600, 6.13)
            };
        
        \node[blue] at (axis cs:200,43) {(5444)};
        \node[blue] at (axis cs:400,37.2) {(8846)};
        \node[blue] at (axis cs:870,11) {(5615)};
        \node[blue] at (axis cs:1610,7) {(6765)};

        \addplot[thick,color=blue,mark=*,dashed] coordinates {
            (200, 33.00)
            (400, 30.20)
            (800, 8.16)
            (1600, 5.02)
            };

        \addplot[thick,color=orange,mark=diamond*] coordinates {
            (200, 35.36)
            (400, 18.15)
            (800, 8.49)
            (1600, 4.52)
            };
            
        \node[orange] at (axis cs:200,25.5) {(3268)};
        \node[orange] at (axis cs:390,12.5) {(4261)};
        \node[orange] at (axis cs:790,5.8) {(4089)};
        \node[orange] at (axis cs:1600,3) {(4108)};

        \addplot[thick,color=orange,mark=diamond*,dashed] coordinates {
            (200, 29.12)
            (400, 14.79)
            (800, 6.61)
            (1600, 3.38)
            };

        \addplot[thick,color=black,dotted,opacity=0.5] coordinates {
            (200, 38.34)
            (400, 19.17)
            (800, 9.58)
            (1600, 4.79)
        };
        \addplot[thick,color=black,dotted,opacity=0.5] coordinates {
            (200, 33.00)
            (400, 16.50)
            (800, 8.25)
            (1600, 4.12)
        };
        \addplot[thick,color=black,dotted,opacity=0.5] coordinates {
            (200, 35.36)
            (400, 17.68)
            (800, 8.84)
            (1600, 4.42)
        };
        \addplot[thick,color=black,dotted,opacity=0.5] coordinates {
            (200, 29.12)
            (400, 14.56)
            (800, 7.28)
            (1600, 3.64)
        };

        \legend{BJ, BJ (linear solve), BGS, BGS (linear solve)};
    \end{loglogaxis}
    \end{tikzpicture}
    \caption{Strong scaling timings for Case III, the flattened jelly roll problem (112M DoFs). The total number of GMRES iterations are given in parentheses. Ideal scaling is given by dotted lines.}
    \label{fig:scaling_jellyroll}
\end{figure}

For Case~III, the strong scaling results for the flattened jelly-roll geometry are shown in \cref{fig:scaling_jellyroll}.  
The highly anisotropic domain and mesh (\SI{1e-4}{\meter} $\times$ \SI{0.05}{\meter} $\times$ \SI{1}{\meter}) pose additional challenges for solver performance.  
The BJ preconditioner reduces the total solve time from approximately \(2300\)~seconds at 200 processors to \(368\)~seconds at 1600 processors, corresponding to a \(6.3\times\) speedup and a parallel efficiency of \(0.79\) relative to ideal \(8\times\) scaling.  
The BGS preconditioner is consistently faster, decreasing from \(2122\)~seconds to \(271\)~seconds over the same processor range, yielding a \(7.8\times\) speedup and an efficiency of \(0.97\).  
The anisotropy of this case degrades the performance of AMG, resulting in a substantially larger number of GMRES iterations: 30--40 per Newton iteration for BGS and 50--80 per Newton iteration for BJ.  
In some instances, a single Newton step requires 100--300 GMRES iterations.  
To assess whether this iteration growth is intrinsic to the block preconditioning strategy or driven by multigrid's poor performance, we performed additional experiments on smaller problems (not shown) in which AMG was replaced by a sparse direct (LU) solver on the block subproblems.  
In these tests, the number of GMRES iterations was reduced to around 10 per Newton iteration, indicating that the increased iteration counts observed here are primarily attributable to the limitations of AMG in the strongly anisotropic regime.  
For the large-scale simulations, however, the memory requirements of LU factorization make this approach impractical.  
Overall, both preconditioners retain strong parallel scalability, with BGS exhibiting lower computational times due to reduced iteration counts.  
Since the total solve time in this case is dominated by the linear solver, which itself scales well, the overall strong scaling appears more favorable than in the other cases, despite the increased number of Krylov iterations.  
These results suggest that multigrid approaches tailored to highly anisotropic geometries could further improve solver robustness and efficiency for cases like jelly-roll cells.

\begin{figure}[htb!]
\centering
\begin{tikzpicture}
    \begin{loglogaxis}[
      log ticks with fixed point,
      scaled ticks=false,
      tick scale binop=\times,
        width=3.75in,
        xlabel={number of processors},
        ylabel={solve time (min)},
        title={Strong Scaling (TPMS)},
        xtick={200,400,800,1600},
        ytick={2.5,5,10,20},
        log basis x=2,
        log basis y=2,
        grid=both,
        grid style={dotted, gray!30},
        legend style={at={(0.98,0.98)}, anchor=north east},
        ]

        \addplot[thick,color=blue,mark=*] coordinates {
            (200, 21.24)
            (400, 10.85)
            (800, 5.55)
            (1600, 4.54)
            }
            node[pos=0.0, below, yshift=-2mm]{(2283)}
            node[pos=0.33, below, yshift=-1.5mm]{(2214)}
            node[pos=0.67, below, yshift=-2mm]{(2054)}
            node[pos=1.0, below]{(3574)}
            ;

        \addplot[thick,color=blue,mark=*,dashed] coordinates {
            (200, 12.35)
            (400, 5.85)
            (800, 2.78)
            (1600, 2.72)
            };

        \addplot[thick,color=orange,mark=diamond*] coordinates {
            (200, 30.55)
            (400, 12.22)
            (800, 6.35)
            (1600, 4.98)
            }
            node[pos=0.0, above]{(1946)}
            node[pos=0.33, above right]{(1873)}
            node[pos=0.67, above right]{(1757)}
            node[pos=1.0, above]{(2962)}
            ;

        \addplot[thick,color=orange,mark=diamond*,dashed] coordinates {
            (200, 19.85)
            (400, 6.78)
            (800, 3.30)
            (1600, 3.14)
            };

        \addplot[thick,color=black,dotted,opacity=0.5] coordinates {
            (200, 21.24)
            (400, 10.62)
            (800, 5.31)
            (1600, 2.65)
        };
        \addplot[thick,color=black,dotted,opacity=0.5] coordinates {
            (200, 12.35)
            (400, 6.17)
            (800, 3.09)
            (1600, 1.54)
        };
        \addplot[thick,color=black,dotted,opacity=0.5] coordinates {
            (200, 30.55)
            (400, 15.28)
            (800, 7.64)
            (1600, 3.82)
        };
        \addplot[thick,color=black,dotted,opacity=0.5] coordinates {
            (200, 19.85)
            (400, 9.93)
            (800, 4.96)
            (1600, 2.48)
        };

        \legend{BJ, BJ (linear solve), BGS, BGS (linear solve)};
    \end{loglogaxis}
    \end{tikzpicture}
    \caption{Strong scaling timings for Case IV, the interdigitated gyroid electrodes problem (237M DoFs). The total number of GMRES iterations are given in parentheses. Ideal scaling is given by dotted lines.}
    \label{fig:scaling_tpms}
\end{figure}
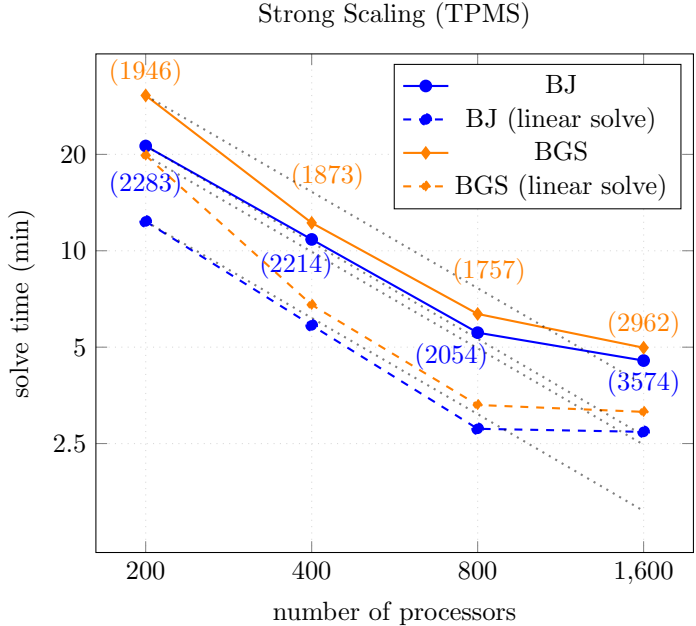

Finally, for Case~IV, the strong scaling results for the interpenetrating gyroid geometry are shown in \cref{fig:scaling_tpms}.  
For the BJ preconditioner, the total solve time decreases from approximately \(1274\)~seconds at 200 processors to \(272\)~seconds at 1600 processors, corresponding to a \(4.7\times\) speedup and a parallel efficiency of \(0.59\).  
The BGS preconditioner remains slower across all processor counts, reducing the solve time from \(1833\)~seconds to \(299\)~seconds over the same range, yielding a \(6.1\times\) speedup and an efficiency of \(0.76\).  
However, both preconditioners show increased sensitivity to mesh partitioning at higher core counts, particularly at 1600 processors where the iteration counts rise sharply (from \(2054\) to \(3574\) for BJ and from \(1757\) to \(2962\) for BGS).  
This sensitivity is consistent with the complex geometry of the TPMS and with the use of a tetrahedral mesh, which results in fewer degrees of freedom per element and therefore exposes communication bottlenecks sooner than in the hexahedral-mesh cases.

\section{Conclusion}
\label{sec:conclusion}

We have developed and evaluated scalable block preconditioning strategies for the fully implicit solution of the pseudo-4D Doyle--Fuller--Newman (DFN) battery model in three-dimensional domains.
By leveraging the block structure of the coupled electrode- and particle-scale equations, we constructed block-diagonal and block-triangular preconditioners that combine algebraic multigrid for electrode-level operators with efficient direct solvers for the particle-scale diffusion blocks.

Extensive scalability experiments were conducted on a range of battery cell geometries, including homogeneous and heterogeneous cubic cells, flattened jelly-roll architectures, and interpenetrating gyroid electrodes.
In all cases, the proposed preconditioners enabled robust convergence of GMRES and strong parallel scalability, with efficient solution of systems containing up to hundreds of millions of unknowns.
The block-diagonal approach usually provided lower wall-clock times, while the block-triangular variant yielded reduced iteration counts, especially for highly anisotropic domains.

These results demonstrate that fully coupled, three-dimensional pseudo-4D DFN simulations can be performed efficiently and robustly at scale using iterative solvers equipped with appropriate block preconditioners.
This advances the practical applicability of high-fidelity battery modeling for realistic cell architectures and operating conditions.

\appendix
\section{Implementation details for the particle equation}
\label{sec:appendix}

Firedrake requires equations to be solved using the finite element method.
Further, there is no obvious way to implement a 4th pseudo-dimension for the particle equation.
Therefore, we will consider each discrete value of the discretized solid-phase concentration as a finite element function, and manually implement a finite difference scheme using these functions.

We define a vector-valued function space \( \mathcal{V}_c = V_h^{N_c} \), where each component corresponds to the solid-phase concentration at one radial location.
The solid-phase concentration field is represented as \( \mathbf{c}_s \in \mathcal{V}_c \), with test functions \( \mathbf{v} \in \mathcal{V}_c \) and components \( c_i \), \( v_i \).

We formulate the weak form by enforcing the finite difference residual at each radial node in a variational sense:
\begin{equation}
\label{eq:particle-weak}
\int_{\Omega_n \cup \Omega_p} \left( \frac{\partial c_i}{\partial t} - \mathcal{L}_h(\mathbf{c}_s)_i \right) v_i \, \mathrm{d}\mathbf{x} = 0, \quad \text{for } i=1,\dots,N_c,
\end{equation}
where \( \mathcal{L}_h(\mathbf{c}_s)_i \) denotes the finite difference approximation of the radial Laplacian, including boundary conditions (see \cref{eq:Lh,eq:L0,eq:LR}).

Therefore, the fully coupled variational problem is: Find \( \mathbf{c}_s \in \mathcal{V}_c \), \( c_e, \phi_e \in V_h \), and \( \phi_s \in V_h^0 \) such that the following \cref{eq:particle-weak,eq:ce-weak,eq:phie-weak,eq:phis-weak} are satisfied for all test functions \( \mathbf{v} \in \mathcal{V}_c \), \( v_{c_e}, v_{\phi_e} \in V_h \), and \( v_{\phi_s} \in V_h^0 \).

The weak form \cref{eq:particle-weak} is implemented directly in Firedrake by defining \( \mathbf{c}_s \) using a \texttt{VectorFunctionSpace}, evaluating the finite difference operators symbolically over the components, and associating one test function per radial node.
Because the number of radial points \( N_c \) is small (typically 5–20), this approach remains computationally efficient and maintains a unified finite element representation of the full battery system.

\section{Solver parameters}
\label{sec:solver_parameters}

Here we provide the PETSc options used to construct the solvers for described in \cref{sec:solvers}.
Note that we are using the default parameters for the nonlinear solver.

\begin{lstlisting}[
    label={lst:dfn_petsc_options},
    language=Python,
    frame=single,
    basicstyle=\footnotesize\ttfamily,
    numbers=left,
    numbersep=5pt,
    caption={PETSc solver options.}
]
% Outer Krylov solver
"ksp_type": "gmres",
"ksp_gmres_restart": 30, # or 100 for Case III

% Block preconditioner
"pc_type": "fieldsplit",
"pc_fieldsplit_type": "multiplicative", % or "additive"
"mat_type": "aij",

% Ordering used for the results in this paper
"fieldsplit_0_fields": i_phie, # index for phi_e
"fieldsplit_1_fields": i_cs, # index for c_s
"fieldsplit_2_fields": i_phis, # index for phi_s
"fieldsplit_3_fields": i_ce, # index for c_e

% AMG parameters. Repeat for each electrode-level block
% For Cases I, II, IV
"fieldsplit_[i]": { % replace i by i_ce, i_phie, and i_phis
    "pc_type": "hypre",
    "pc_hypre_boomeramg": {
        "strong_threshold": 0.7,
        "coarsen_type": "HMIS",
        "agg_nl": 3,
        "interp_type": "ext+i",
        "agg_num_paths": 5,
    },
},
% For Case, III
"fieldsplit_[i]": {
    "pc_type": "hypre",
        "pc_hypre_boomeramg": {
          "strong_threshold": 0.9,
          "coarsen_type": "HMIS",
          "agg_nl": 0,
          "interp_type": "ext+i",
          },
      }


% Preconditioner for particle system
"fieldsplit_[i_cs]": {
    "pc_type": "pbjacobi",
},
\end{lstlisting}

\section{Model Parameters}
\label{sec:model_parameters}

This appendix summarizes all parameters used in the Doyle--Fuller--Newman (DFN) model described in \cref{sec:gov}.
All parameter values are from \cite{marquis2019asymptotic}, taken directly from PyBaMM \cite{Sulzer2021}.

\begin{table}[h!]
\centering
\caption{Universal physical constants.}
\label{tab:constants}
\begin{tabular}{lll l}
\hline
Symbol & Description & Value & Unit \\ \hline
$F$ & Faraday constant & 96485.33 & \si{\coulomb\per\mole} \\
$R$ & Ideal gas constant & 8.31446 & \si{\joule\per\mole\per\kelvin} \\
$T$ & Temperature & 298.15 & \si{\kelvin} \\
\hline
\end{tabular}
\end{table}

\begin{table}[h!]
\centering
\caption{Electrolyte parameters.}
\label{tab:electrolyte}
\begin{tabular}{lll l}
\hline
Symbol & Description & Value & Unit \\ \hline
$c_{e,0}$      & Initial electrolyte concentration & 1000  & \si{\mole\per\cubic\meter} \\
$t_+^0$        & Cation transference number        & 0.4   & - \\
$b$            & Bruggeman exponent                & 1.5   & - \\
$D_e(c_e,T)$   & Bulk electrolyte diffusivity      & see \cref{eq:De_app} & - \\
$\kappa(c_e,T)$& Bulk electrolyte conductivity     & see \cref{eq:kappa_app} & - \\
\hline
\end{tabular}
\end{table}

\begin{align}
D_e(c_e, T)
&= 5.34 \times 10^{-10}
\exp\!\left(-0.65 \frac{c_e}{1000}\right)
\exp\!\left(\frac{E_{D_e}}{R}\left(\frac{1}{298.15} - \frac{1}{T}\right)\right),
\label{eq:De_app} \\[0.5em]
\kappa(c_e, T)
&= \exp\!\left(\frac{34700}{R}\left(\frac{1}{298.15} - \frac{1}{T}\right)\right)
\left(0.0911 + 1.9101 c - 1.052 c^2 + 0.1554 c^3\right),
\label{eq:kappa_app}
\end{align}
where $c = c_e / 1000$ is the electrolyte concentration in \si{\mole\per\liter} and \linebreak
$E_{D_e} = \SI{37040}{\joule\per\mole}$.

\begin{table}[h!]
\centering
\caption{Cathode (positive electrode) parameters.}
\label{tab:cathode}
\begin{tabular}{lll l}
\hline
Symbol & Description & Value & Unit \\ \hline
$\varepsilon_s$         & Active material volume fraction      & \num{0.5}          & - \\
$\varepsilon$           & Porosity                            & \num{0.3}          & - \\
$\sigma$                & Solid-phase conductivity            & \num{10}           & \si{\siemens\per\meter} \\
$k$                     & Reaction rate constant              & \num{6e-7}         & \si{\ampere\per\square\meter\,(\meter^3\per\mole)^{1.5}} \\
$\alpha_a = \alpha_c$   & Charge-transfer coefficients        & \num{0.5}          & - \\
$c_{s,\max}$            & Maximum solid-phase concentration   & \num{5.12e4}       & \si{\mole\per\cubic\meter} \\
$c_{s,0}$               & Initial solid-phase concentration   & \num{3.07e4}       & \si{\mole\per\cubic\meter} \\
$a$                     & Specific interfacial area           & \num{1.5e5}        & \si{\square\meter\per\cubic\meter} \\
$D_s$                   & Solid diffusivity                   & \num{1.0e-13}      & \si{\square\meter\per\second} \\
$R_s$                   & Particle radius                     & \num{1.0e-5}       & \si{\meter} \\
$U_\mathrm{ocp}$        & Open-circuit potential              & see \cref{eq:ocp_cathode} & - \\
\hline
\end{tabular}
\end{table}

\begin{table}[h!]
\centering
\caption{Anode (negative electrode) parameters.}
\label{tab:anode}
\begin{tabular}{lll l}
\hline
Symbol & Description & Value & Unit \\ \hline
$\varepsilon_s$         & Active material volume fraction      & \num{0.6}          & - \\
$\varepsilon$           & Porosity                            & \num{0.3}          & - \\
$\sigma$                & Solid-phase conductivity            & \num{100}          & \si{\siemens\per\meter} \\
$k$                     & Reaction rate constant              & \num{2e-5}         & \si{\ampere\per\square\meter\,(\meter^3\per\mole)^{1.5}} \\
$\alpha_a = \alpha_c$   & Charge-transfer coefficients        & \num{0.5}          & - \\
$c_{s,\max}$            & Maximum solid-phase concentration   & \num{2.50e4}       & \si{\mole\per\cubic\meter} \\
$c_{s,0}$               & Initial solid-phase concentration   & \num{2.00e4}       & \si{\mole\per\cubic\meter} \\
$a$                     & Specific interfacial area           & \num{1.8e5}        & \si{\square\meter\per\cubic\meter} \\
$D_s$                   & Solid diffusivity                   & \num{3.9e-14}      & \si{\square\meter\per\second} \\
$R_s$                   & Particle radius                     & \num{1.0e-5}       & \si{\meter} \\
$U_\mathrm{ocp}$        & Open-circuit potential              & see \cref{eq:ocp_anode} & - \\
\hline
\end{tabular}
\end{table}

\begin{table}[h!]
\centering
\begin{threeparttable}
\caption{Separator parameters.}
\label{tab:separator}
\begin{tabular}{lll l}
\hline
Symbol & Description & Value & Unit \\ \hline
$\varepsilon$  & Porosity                 & \num{1.0}        & - \\
$\sigma$       & Solid-phase conductivity & \num{0}$^\dagger$ & \si{\siemens\per\meter} \\
\hline
\end{tabular}
\begin{tablenotes}
\small
\item[$\dagger$] This is the physical value; a small non-zero value is chosen for the effective solid-phase conductivity $\sigma^\mathrm{eff}$ for numerical reasons. See \cref{sec:implementation} for details.
\end{tablenotes}
\end{threeparttable}
\end{table}

Effective transport coefficients are obtained via Bruggeman corrections as described in \cref{sec:gov}.

The interfacial kinetics use $U_{\mathrm{ocp}}(c_s^{\mathrm{surf}})$ evaluated as a function of the particle surface concentration.
Let the (surface) stoichiometry
\[
\theta = \frac{c_s^{\mathrm{surf}}}{c_{s,\max}}.
\]

The cathode (LiCoO\textsubscript{2}) open-circuit potential is
\begin{align}
U_{\mathrm{ocp}}(\theta)
&=
2.16216
+ 0.07645 \tanh\!\left(30.834 - 54.4806\,\tilde{\theta}\right)
\nonumber\\
&\quad
+ 2.1581 \tanh\!\left(52.294 - 50.294\,\tilde{\theta}\right)
\nonumber\\
&\quad
- 0.14169 \tanh\!\left(11.0923 - 19.8543\,\tilde{\theta}\right)
+ 0.2051 \tanh\!\left(1.4684 - 5.4888\,\tilde{\theta}\right)
\nonumber\\
&\quad
+ 0.2531 \tanh\!\left(\frac{-\tilde{\theta} + 0.56478}{0.1316}\right)
- 0.02167 \tanh\!\left(\frac{\tilde{\theta} - 0.525}{0.006}\right),
\label{eq:ocp_cathode}
\end{align}
where $\tilde{\theta} = 1.062\,\theta$.

The anode (graphite) open-circuit potential is
\begin{align}
U_{\mathrm{ocp}}(\theta)
&=
0.194
+ 1.5 \exp\!\left(-120\,\theta\right)
+ 0.0351 \tanh\!\left(\frac{\theta - 0.286}{0.083}\right)
\nonumber\\
&\quad
- 0.0045 \tanh\!\left(\frac{\theta - 0.849}{0.119}\right)
- 0.035 \tanh\!\left(\frac{\theta - 0.9233}{0.05}\right)
\nonumber\\
&\quad
- 0.0147 \tanh\!\left(\frac{\theta - 0.5}{0.034}\right)
- 0.102 \tanh\!\left(\frac{\theta - 0.194}{0.142}\right)
\nonumber\\
&\quad
- 0.022 \tanh\!\left(\frac{\theta - 0.9}{0.0164}\right)
- 0.011 \tanh\!\left(\frac{\theta - 0.124}{0.0226}\right)
\nonumber\\
&\quad
+ 0.0155 \tanh\!\left(\frac{\theta - 0.105}{0.029}\right).
\label{eq:ocp_anode}
\end{align}

\begin{table}[h!]
\centering
\caption{Current collector parameters for Case III.}
\label{tab:current_collectors}
\begin{tabular}{lll l}
\hline
Material & Description & Value & Unit \\ \hline
Aluminum & Electrical conductivity & \num{3.77e7} & \si{\siemens\per\meter} \\
Copper   & Electrical conductivity & \num{5.96e7} & \si{\siemens\per\meter} \\
\hline
\end{tabular}
\end{table}

To define the applied current corresponding to a 1C discharge, we compute the \emph{theoretical} (i.e., maximum) cell capacity from the active-material lithium inventory in each porous electrode.
Let $V_n$ and $V_p$ be the volume of the negative and positive electrode domains, respectively.
Given the solid volume fraction $\varepsilon_{s,i}$ and the maximum lithium concentration in the solid phase $c_{s,i}^{\max}$ (in $\mathrm{mol\,m^{-3}}$ of solid), the theoretical capacity of electrode $i$ (in ampere-hours) is computed as
\begin{equation}
Q_i \;=\; \frac{F}{3600}\,\varepsilon_{s,i}\,V_i\,c_{s,i}^{\max},
\qquad i\in\{n,p\},
\label{eq:theoretical-capacity-electrode}
\end{equation}
where constant and the factor $3600$ converts coulombs to ampere-hours.
The theoretical cell capacity is taken as the minimum of the electrode capacities,
\begin{equation}
Q \;=\; \min\{Q_n,Q_p\},
\label{eq:theoretical-capacity-cell}
\end{equation}
so that the limiting electrode determines the maximum cyclable lithium inventory.

Finally, for a prescribed C-rate $\mathsf{C}$, we define the applied current magnitude as
\begin{equation}
I_{\mathrm{app}} \;=\; \mathsf{C}\,Q,
\label{eq:applied-current-crate}
\end{equation}
where $\mathsf{C}=1$ is chosen for all simulations.

Note that in our configuration, the electrodes are assigned equal geometric volumes and solid volume fractions, which results in $Q_n < Q_p$ due to the smaller maximum lithium concentration of the graphite anode.
In contrast, commercial lithium-ion cells typically employ an oversized anode, so that the cathode is capacity-limiting and $Q_p < Q_n$.


\bibliographystyle{siamplain}
\bibliography{references}
\end{document}